\documentclass[12pt,amsart]{article}
\usepackage{amsmath, amssymb, amsfonts, amsthm,amscd,xypic}
\newtheorem{theorem}{Theorem}
\newtheorem{lemma}[theorem]{Lemma}

\newtheorem{proposition}[theorem]{Proposition}
\newtheorem{corollary}[theorem]{Corollary}

\theoremstyle{definition}
\newtheorem{remark}[theorem]{Remark}

\newtheorem{definition}[theorem]{Definition}
\newtheorem{question}[theorem]{Question}

 \newcommand{\SL}{\operatorname{SL}}

  \newcommand{\al}{\alpha}
\newcommand{\defpar}{\tau}
    \newcommand{\Del}{\Delta}
     
     \newcommand{\Lam}{\Lambda}
  
   \newcommand{\wed}{\wedge}

 \newcommand{\ip}[1]{\langle #1 \rangle}
 \newcommand{\haq}{\hat{Q}}
  \renewcommand{\tilde}{\widetilde}

  \def\b1{\text{\bf\large 1}}  

  \def\ip<#1>{\langle#1\rangle}   

  \newcommand{\codim}{\operatorname{codim}}

\newcommand{\gr}{\operatorname{gr}}
  \newcommand{\Hom}{\operatorname{Hom}}
  \newcommand{\Imo}{\operatorname{Im}}
\newcommand{\Gr}{\operatorname{Gr}}

\newcommand{\Ad}{\operatorname{Ad}}

  \newcommand\elal{\Bbb{L}}
    \newcommand\xx{\Bbb{X}}
    \newcommand\exx{\Bbb{X}}
    \newcommand\gm{\Bbb{G}_m}

\newcommand{\beqn}{\begin{equation}}
\newcommand{\eeqn}{\end{equation}}

\newcommand{\tset}[1]{\underset{#1}{\times}}

 \newcommand{\fb}{\mathfrak{b}}
 \newcommand{\fg}{\mathfrak{g}}
 \newcommand{\fh}{\mathfrak{h}}
 \newcommand{\fl}{\mathfrak{l}}
 
 \newcommand{\fp}{\mathfrak{p}}
 \newcommand{\fs}{\mathfrak{s}}
 \newcommand{\ft}{\mathfrak{t}}
 \newcommand{\fu}{\mathfrak{u}}

\def\u.{^{\bullet}}

\newcommand{\bc}{\mathbb{C}}

\newcommand{\bz}{\mathbb{Z}}

\newcommand{\frg}{\mathfrak{g}}
\newcommand{\frh}{\mathfrak{h}}

\newcommand{\frk}{\mathfrak{k}}
\newcommand{\frl}{\mathfrak{l}}

\newcommand{\frp}{\mathfrak{p}}

\newcommand{\fru}{\mathfrak{u}}

\newcommand\mt{\mathcal{T}}

 \newcommand{\ca}{\mathcal{A}}

 \newcommand{\cf}{\mathcal{F}}

 \newcommand{\cl}{\mathcal{L}}

\begin{document}

  \title{Eigenvalue problem and a new product in cohomology of flag
varieties}
  \author{Prakash Belkale and Shrawan Kumar\\Department of Mathematics\\
University of North Carolina\\
Chapel Hill, NC  27599--3250}
  \maketitle
 \section{Introduction}
Let $G$ be a connected semisimple complex algebraic group and let $P$ be a parabolic subgroup. In this paper we define a new (commutative and associative) product $\odot_0$ on the cohomology of the homogenous spaces $G/P$ and use this to give a more efficient solution of the eigenvalue problem and also for the problem of determining the existence  of
$G$-invariants in the tensor product of irreducible representations of $G$. On the other hand, we show that this new product is intimately connected with the Lie algebra cohomology of the nil-radical of $P$ via some works of Kostant and Kumar. We also initiate a uniform study of the geometric Horn problem for an arbitrary group $G$ by obtaining two (a priori) different sets of necessary recursive conditions  to determine  when a cohomology product of Schubert classes in
 $G/P$ is non-zero. Hitherto, this was studied largely only for the group $\SL(n)$ and its maximal parabolics $P$.

This new cohomology product is a certain deformation of the classical product. If $w\in W^P$ ($W^P$ is the set of minimal length representatives in the cosets of $W/W_P$), let $[\bar{\Lambda}_w^P]\in H^*(G/P)$ be the cohomology class of the subvariety $\bar{\Lambda}^P_w:=\overline{w^{-1}BwP}\subseteq G/P$. If the structure coefficients for the classical product are written as
$$[\bar{\Lambda}^P_u] \cdot [\bar{\Lambda}^P_v]=\sum_w c^w_{u,v}[\bar{\Lambda}^P_w],$$
then the new product $\odot_0$ is a restricted sum
$$[\bar{\Lambda}^P_u] \odot_0 [\bar{\Lambda}^P_v]={\sum_w}' c^w_{u,v}[\bar{\Lambda}^P_w],$$ where the summation is only over a smaller set of $w$ which
satisfy a certain numerical condition involving $u$, $v$ and $w$
(Definition  \ref{levicohomology}). This numerical
condition is best understood in terms of  a
notion which we call $L$-movability (c.f., Definition \ref{D1} and Section \ref{ref11}). If $P$ is a minuscule maximal
parabolic (e.g., for any maximal parabolic in $\SL(n)$), $\odot_0$ coincides with the 
classical product. However, even for any nonmaximal 
parabolic $P$ in $\SL(n)$, $\odot_0$ differs from the classical product. 
\subsection{Eigenvalue problem}
Choose a Borel subgroup $B$ and a maximal torus $H\subset B$ of $G$. Let $K$ be a maximal compact subgroup of $G$ chosen such that
 $i\frh_{\Bbb{R}}$ is the Lie algebra of
a maximal torus of $K$, where $\frh_{\Bbb{R}}$ is a real form of
the Lie algebra $\fh$ of $H$. There is a natural homeomorphism
$C:\frk/K\to \frh_{+},$ where $K$ acts on $\frk$ by the adjoint
representation and $\frh_{+}$ is the positive Weyl chamber in
$\frh_{\Bbb{R}}$.

The {\it {eigenvalue problem}} is concerned with the following question:
\begin{enumerate}
\item[(E)]  Determine all the  $s$-tuples $(h_1,\dots,h_s)\in \frh_{+}^s$ for which there exist $(k_1,\dots,k_s)\in\frk^s$ such that  $C(k_j)=h_j$ for
$j=1,\dots,s$, and
$$\sum_{j=1}^s k_j=0.$$
\end{enumerate}
Given a (standard) maximal parabolic subgroup $P$, let $\omega_P$ denote the corresponding fundamental weight. This is invariant under the Weyl group of $P$.

For (E), Leeb-Millson  (following the works of Klyachko [Kl], Belkale
[Bel1] and Berenstein-Sjamaar ~\cite{BerensteinSjamaar}) obtained the following:

{\it Let $(h_1, \cdots, h_s)\in\frh_{+}^s$. Then, the following are equivalent:}
\begin{enumerate} \item[(A)] $\exists (k_1, \cdots, k_s)\in \frk^s$ {\it such 
that}
$\,k_1+ \cdots+ k_s=0,$ {\it and }\, $C(k_j)=h_j,$\, {\it for } $ j=1, \cdots, s.$
\item[(B)] {\it For every standard maximal parabolic subgroup $P$ in $G$ and every choice 
of
$s$-tuples $(w_1, \cdots, w_s)\in (W^P)^s$ such that} $$[{\bar{\Lambda}}^P_{w_1}]\cdot
\,\cdots \,\cdot [{\bar{\Lambda}}^P_{w_s}] = [\bar{\Lambda}^P_{e}] \in H^*(G/P, \Bbb
Z),$$ \end{enumerate}
{\it the following inequality holds:}
 \begin{equation}\label{IntroIneq}
\omega_P(\sum_{j=1}^s\,w_j^{-1}h_j)\leq 0. \end{equation} 

In our Theorem ~\ref{eigen}, we show that we can replace the condition (B) by a smaller set of inequalities. This is one of the main theorems of this paper.

\noindent
{\bf Theorem:} (A) {\it is equivalent to}
\begin{enumerate}
\item[(B$'$)] {\it For every standard maximal parabolic subgroup $P$ in $G$ and every 
choice of $s$-tuples  $(w_1, \cdots, w_s)\in (W^P)^s$ such that}
$$[\bar{\Lambda}^P_{w_1}]\odot_0 \,\cdots \,\odot_0 [\bar{\Lambda}^P_{w_s}] = 
[\bar{\Lambda}^P_{e}] \in (H^*(G/P, \Bbb Z),\odot_0),$$
\end{enumerate}
{\it the inequality (\ref{IntroIneq}) holds}.

 Knutson-Tao-Woodward [KTW] proved that for $G=\SL(n)$ the inequalities (B) (which are,
in this case, the same set of inequalities as (B$'$) because of Lemma \ref{minuscule})
 are irredundant. However, in the case of the groups of type $B_3$ and $C_3$,
Kumar-Leeb-Millson ~\cite{KumarLeebMillson} proved that of the total of 135 inequalities
corresponding to the system (B), 33 are redundant. In contrast, making use of the 
tables in Section 10, it can be seen that none of the inequalities
given by (B$'$) are redundant for either $B_3$ or $C_3$.

Let $\nu_1, \cdots, \nu_s$ be dominant integral weights. Then, consider the problem of finding conditions on $\nu_j$'s such that the space of $G$-invariants
$$H^0((G/B)^s,  \cl(N\nu_1)\boxtimes \cdots \boxtimes\cl(N\nu_s))^G$$ is nonzero for some $N>0$. As is well known, by using symplectic geometry, this
problem is equivalent to the eigenvalue problem. For a precise solution of this problem, see Theorem \ref{EVT}.

\subsection{Relation to Lie algebra cohomology}
Using a celebrated theorem of Kostant  on the cohomology of some nilpotent Lie algebras (Theorem ~\ref{Th9.A}), we relate the new cohomology product
$\odot_0$ to Lie algebra cohomology. More specifically, in Theorem ~\ref{liecohomology}, we exhibit an explicit isomorphism of graded rings
 $$\phi : \bigl( H^*(G/P, \bc ), \odot_0\bigr) \simeq \bigl[
H^*(\fu_P)\otimes H^*(\fu_P^-)\bigr]^{\fl},$$
where we take the tensor product algebra structure on
the right side. This isomorphism is used  to determine the
structure of $(H^*(G/B),\odot_0)$ completely (cf., Corollary ~\ref{leviprod}).

\subsection{Geometric Horn problem} 
One of the consequences of the work of Klyachko [Kl],
and the saturation theorem of Knutson-Tao ~\cite{KT} is that one can tell when a product
of Schubert classes in a given Grassmannian $\Gr(r,n)$ (which is a homogenous space for
$\SL(n)$ corresponding to a maximal parabolic) is nonzero by writing down a series of
inequalities coming from knowing the answer to the same question for smaller
Grassmannians. These inequalities are the eigenvalue inequalities (\ref{IntroIneq}) for 
$\SL(r)$ where
the conjugacy classes are determined from the Schubert classes. We refer the reader to
the survey article of Fulton [Fu2] for details.

In this paper we initiate a uniform study of solving the same problem for arbitrary semisimple groups.

We obtain two (a priori) different sets of necessary conditions (Theorems ~\ref{T2} and
~\ref{weakHorn}) for a product (with any number of factors) of Schubert classes to be
nonzero in the cohomology of $G/P$. The inequalities given by Theorem ~\ref{T2}, which we
call {\it character inequalities}, are similar to the eigenvalue inequalities for the
Levi subgroup $L$ of $P$. These inequalities provide the first known numerical criterion
for the vanishing of intersection multiplicities beyond the codimension condition for an
arbitrary semisimple $G$.

The second set of inequalities, which we call the {\it dimension inequalities}, given by
Theorem ~\ref{weakHorn} are based on dimension counts. In the case of $G=\SL(n)$ and $P$
a maximal parabolic, the dimension and character inequalities coincide. In general, even
for minuscule maximal parabolics, these inequalities are (a priori) different.

The character inequalities can be refined further if we consider the new product
$\odot_0$ instead of the classical product (Theorem ~\ref{T2'}). The character
inequalities for $\odot_0$ are sufficient for the homogenous spaces of the form $G/B$ by
virtue of Corollary ~\ref{leviprod}. \subsection{Examples} In Section ~\ref{exemples} we
give the multiplication tables under the deformed product $\odot$ for $G/P$ for all the
rank $3$ complex simple groups $G$ and maximal parabolic subgroups $P$.

\vskip2ex

Sections 7,8 and 9 are independent of each other.
\vskip4ex
\noindent
{\bf Acknowledgements:}  Both the 
authors thank the NSF for partial supports. We are
pleased to thank
 A. Knutson for a helpful conversation leading to Theorem \ref{liecohomology}.

\section{Notation and Preliminaries on the Intersection Theory in $G/P$}\label{introo}

Let $G$ be a connected reductive complex algebraic group. We choose a Borel subgroup $B$ and a maximal torus $H\subset B$  and let $W=W_G:=N_G(H)/H$ be the associated Weyl group, where $N_G(H)$ is the normalizer of $H$ in $G$.  Let $P\supseteq B$ be a standard parabolic subgroup of $G$ and let $U=U_P$ be its unipotent radical. Consider the  Levi subgroup $L=L_P$ of $P$ containing $H$, so that $P$ is the semi-direct product of $U$ and $L$. Then,  $B_L:=B\cap L$ is a Borel subgroup of
$L$. Let $X(H)$ denote the character group of $H$, i.e., the group of all the algebraic group morphisms $H \to \Bbb G_m$. Then, $B_L$ being the semidirect product of its commutator $[B_L,B_L]$ and $H$, any $\lambda \in X(H)$ extends uniquely to a character of $B_L$. Similarly, for any algebraic subgroup $S$ of $G$, let $O(S)$ denote the set of all the one parameter subgroups in $S$, i.e., algebraic group morphisms $\Bbb G_m \to S$. We denote the Lie algebras of $G,B,H,P,U,L,B_L$ by the corresponding Gothic characters: $\fg,\fb,\fh,\fp,\fu,\fl,\fb_L$ respectively.  Let $R=R_\fg$ be
the set of roots of $\fg$ with respect to the Cartan subalgebra $\fh$ and let $R^+$ be the set of positive roots (i.e., the set of roots of $\fb$).
Similarly, let $R_\fl$ be the set of roots of $\fl$ with respect to $\fh$ and
$R_\fl^+$ be the set of roots of $\fb_L$. Let  $\Delta = \{\alpha_1, \cdots, \alpha_\ell\} \subset R^+$ be the set of simple
roots, where $\ell$ is the semisimple rank of $G$ (i.e., the dimension of
$\fh':=\fh\cap [\fg,\fg]$). We denote by $\Delta(P)$ the set of simple roots contained in $R_\fl$.  For any $ 1\leq
j\leq \ell$, define the element $x_j\in \fh'$
by
\begin{equation}\label{eqn0}\alpha_i(x_{j})=\delta_{i,j},\text{ }\forall\text{ } 1\leq i\leq \ell.
\end{equation}

For an $H$-invariant subspace $V\subset\frg$ under the adjoint action, let $R(V)$ denote the set of roots
of $\frg$  appearing in $V$, i.e., $V=\oplus_{\alpha\in R(V)}\,\frg_\alpha$.

Recall that if $W_P$ is the Weyl group of $P$ (which is, by definition, the  Weyl Group of $L$), then in each coset of $W/W_P$ we have a unique member $w$ of minimal length. This satisfies (cf. [Ku2, Exercise 1.3.E]):
\begin{equation}\label{eqn1}
wB_L w^{-1} \subseteq B.
\end{equation}
Let $W^P$ be the set of the minimal length representatives
in the cosets of $W/W_P$.

For any $w\in W^P$, define the (shifted) Schubert cell:
  \[
\Lambda^P_w := w^{-1} BwP \subset G/P.
  \]
Then, it is a locally closed subvariety of $G/P$ isomorphic with the affine
space $\Bbb A^{\ell(w)}, \ell(w)$ being the length of $w$ (cf. [J, Part II,
$\S$13.1]). Its closure is denoted by $\bar{\Lambda}^P_w$, which is an irreducible (projective) subvariety
of $G/P$ of dimension $\ell(w)$. Considered as an element in the group of rational equivalence classes $A_{\ell(w)}(G/P)$ of algebraic degree ${\ell(w)}$ on
$G/P$, it is denoted by  $[\bar{\Lambda}^P_w]$. Since $G/P$ is smooth, setting
$A^*(G/P):= A_{\dim G/P-*}(G/P)$, $A^*(G/P)$ is a commutative associative graded ring under the intersection product $\cdot$ (cf. [Fu1, $\S$8.3]).

Let
$\mu(\bar{\Lambda}^P_w)$ denote the fundamental class of $\bar{\Lambda}^P_w$ considered as an element of the singular homology with integral coefficients $H_{\ell(w)}(G/P, \Bbb Z)$ of $G/P$. Then, from the Bruhat decomposition, the elements
$\{\mu(\bar{\Lambda}^P_w)\}_{w\in W^P}$ form a $\Bbb Z$-basis of
 $H_*(G/P, \Bbb Z)$. Let $\{\epsilon^P_w\}_{w\in W^P}$ be the dual basis of the singular
cohomology with integral coefficients $ H^*(G/P, \Bbb Z)$, i.e., for any $v,w\in W^P$ we have
\[\epsilon^P_v(\mu(\bar{\Lambda}^P_w))=\delta_{v,w}.\]
By [Fu1, Example 19.1.11(b)] and [KLM, Proposition 2.6], there is a graded ring isomorphism, the cycle class map,
\[c:A^*(G/P) \to H^*(G/P, \Bbb Z), \,\,\, [\bar{\Lambda}^P_w]\mapsto
\epsilon^P_{w_oww_o^P},\]
where $w_o$ (resp. $w_o^P$) is the longest element of $W$ (resp. $W_P$).
(For $w\in W^P$, we have $w_oww_o^P\in W^P$ by [KLM, Proposition 2.6].)
{\it From now on, we will identify $A^*(G/P)$ with  $H^*(G/P, \Bbb Z)$ under $c$.} Thus, under the identification $c$,  $ H^*(G/P, \Bbb Z)$ has two $\Bbb Z$-bases: $\{[\bar{\Lambda}^P_w]=\epsilon^P_{w_oww_o^P}\}_{w\in W^P}$ and
$\{\epsilon^P_{w}\}_{w\in W^P}$.

Let $T^P=T(G/P)_e$ be the tangent space of $G/P$ at $e\in G/P$. It carries a
canonical  action of $P$.  For $w\in W^P$, define $T_w^P$ to be the tangent space of $\Lambda_w^P$ at $e$. We shall abbreviate $T^P$ and $T_w^P$ by $T$ and $T_w$ respectively when the reference to $P$ is clear. By (\ref{eqn1}), $B_L$ stabilizes $\Lambda^P_w$ keeping $e$
fixed. Thus,
  \begin{equation}\label{eqn2}  B_L T_w \subset T_w.
\end{equation}

  \begin{lemma}\label{lemma1}
$g \Lambda^P_w$ passes through $e\Leftrightarrow
g\Lambda^P_w = p\Lambda^P_w$ for some $p\in P$.
  \end{lemma}
\begin{proof}  Since $g \Lambda^P_w$ passes through $e$, $g^{-1}\in
 \Lambda^P_w$, i.e., $g\in Pw^{-1}Bw$. Write $g=pw^{-1}bw$, for some $p\in P$ and
$b\in B$. Then, $g \Lambda^P_w = pw^{-1}bw\Lambda^P_w = p\Lambda^P_w $.
  \end{proof}
The following result is the starting point of our analysis.
  \begin{proposition}\label{First}
Take any $s\geq 1$ and any  $(w_1, \cdots, w_s)\in (W^P)^s$  such that
\begin{equation}\label{expected}
\sum_{j=1}^s\,\codim \Lambda_{w_j}^P
 \leq \dim G/P.
 \end{equation}
(Clearly, (3) is equivalent to the following equation:
\begin{equation}\label{expected'}
\sum_{j=1}^s\, \ell(w_j)
 \geq (s-1) \dim G/P.)
 \end{equation}
Then the following three conditions are equivalent:

\begin{enumerate}
\item[(a)] $[\bar{\Lambda}^P_{w_1}]\cdot \,\cdots \,\cdot [\bar{\Lambda}^P_{w_s}] \neq 0 \in A^*(G/P)$.

(Observe that, by the above isomorphism $c$, this is equivalent to the condition:
$\epsilon^P_{w_ow_1w_o^P}\cdots \epsilon^P_{w_ow_sw_o^P}\neq 0.$)
\item[(b)] For generic $(p_1, \cdots, p_s)\in P^s$, the intersection
 \[p_1\Lambda^P_{w_1}\cap \cdots \cap p_s\Lambda^P_{w_s}\]
 is transverse at $e$.
\item[(c)] For generic $(p_1,\dots,p_s)\in P^s$,
 \[\dim(p_1T_{w_1} \cap \cdots \cap p_sT_{w_s}) = \dim G/P -\sum_{j=1}^s\,\codim \Lambda_{w_j}^P .\]
\end{enumerate}
 As proved below, the set of $s$-tuples in (b) as well as (c) is an open subset of $P^s$.
 \end{proposition}
\begin{proof} For $p\in P$ and $w\in W^P$, the tangent space to $p\Lambda^P_w$ at $e\in G/P$ is $pT_w$. Therefore, by the definition of transversality (cf.
[S, Chap. II, $\S$2.1]), (c) is equivalent to  (b). The  subset in (c) is the set of points $(p_1,\dots,p_s)\in P^s$ for which the canonical  morphism
$$T\to\bigoplus_{j=1}^s \frac{T}{p_jT_{w_j}}$$ is  surjective. Therefore, this set is open in $P^s$.

If $(p_1,\dots,p_s)\in P^s$ satisfies the property in (b), the smooth
varieties
 $p_1\Lambda^P_{w_1}, \dots, p_s\Lambda^P_{w_s}$ meet transversally and hence properly at $e\in G/P$. By [Fu1, Proposition 7.1 and Section 12.2]  this implies that (a) holds.

To show that (a) implies (b),  find $g_j\in G$ for $j=1,\dots,s$  so that
$g_1\Lambda^P_{w_1},\dots,g_s \Lambda^P_{w_s}$ meet transversally at a nonempty set of points (cf. Proposition ~\ref{kleiman}). By translation, assume that $e\in G/P$ is one of these points. By Lemma ~\ref{lemma1}, for any $j=1, \cdots, s$, $g_j\Lambda^P_{w_j}=p_j
\Lambda^P_{w_j}$, for some $p_j\in P$. Thus,  we have found a point $(p_1, \dots, p_s)\in P^s$ satisfying the condition in (b). Since the condition is an open condition, we are assured of a\ nonempty open subset of $P^s$ satisfying the property in (b).
\end{proof}

\begin{proposition}[Kleiman]\label{kleiman}
Let a connected algebraic group $G$ act transitively on a smooth
variety $X$ and let  $X_1$, $\dots, X_s$ be irreducible subvarieties of $X$.
Then, there exists a non empty open subset $U\subseteq G^s$ such that for
$(g_1,\dots, g_s)\in U$, the intersection $\bigcap_{j=1}^s \,g_j X_j$ is
proper (possibly empty) and dense in $\bigcap_{j=1}^s \,g_j \bar{X}_j$.

Moreover, if  $X_j$, $j=1,\dots, s$,  are
smooth varieties, we can find such a
$U$ with the additional property that for $(g_1,\dots, g_s)\in U$,
$\bigcap_{j=1}^s \,g_j X_j$ is transverse at each
point of intersection.
\end{proposition}
\begin{proof}
We include a proof of the density statement (the rest of the conclusion is
standard, see [Kle]).

Let  $Y_j:=\bar{X}_j\setminus X_j$. Let $U$ be a
nonempty open subset of $G^s$  such that for $(g_1,\dots, g_s)\in U$, the
following intersections are proper:
\begin{enumerate}
\item $\cap_{j=1}^s \,g_j \bar{X}_j$.
\item For $\ell\in \{1,\dots, s\}$, $\{\cap_{j\in\{1,\dots,
s\}\setminus \{\ell\}} \,\,g_j\bar{X}_j\}\cap g_{\ell} Y_{\ell}.$
\end{enumerate}
 For $(g_1,\dots, g_s)\in U$ and $\ell \in \{1,\dots,s\}$, each irreducible
component of the intersection
$\{\cap_{j\in\{1,\dots, s\}\setminus \{\ell\}} \,\,g_j\bar{X}_j\}\cap
g_{\ell} Y_{\ell}$
is therefore of dimension strictly less than that of  each irreducible
component of $\bigcap_{j=1}^s \,g_j \bar{X}_j$
(since $\dim(Y_{\ell})<\dim(X_{\ell})$). This proves the
density statement.
\end{proof}
  \smallskip

\section{Preliminary Analysis of Levi-movability}\label{Sprelim}

We begin by introducing the central concept of this paper:
\begin{definition}\label{D1}
Let $w_1, \cdots, w_s\in W^P$ be such that
  \begin{equation}\label{dim0}
\sum_{j=1}^s\codim \Lam^P_{w_j}  = \dim G/P.
  \end{equation}
This is equivalent to the condition:
\begin{equation}\label{dim0'}
\sum_{j=1}^s\, \ell(w_j)
 = (s-1) \dim G/P.
 \end{equation}
We then call the $s$-tuple $(w_1, \cdots, w_s)$ {\it Levi-movable} for short $L$-{\it movable} if, for generic $(l_1, \cdots, l_s)\in L^s$, the intersection
$l_1\Lambda_{w_1}\cap \cdots \cap l_s\Lambda_{w_s}$ is transverse at $e$.

Observe that, even though in the definition of $L$-movability we took minimal coset representatives in $W/W_P$, our definition is independent of the choice of coset representatives. Thus, we have the notion of $L$-movability for any $s$-tuple $(w_1, \cdots, w_s)\in (W/W_P)^s.$
\end{definition}

By Proposition ~\ref{First}, if $(w_1, \cdots, w_s)$ is $L$-movable,
 then
$[\bar{\Lambda}^P_{w_1}]\cdot \, \cdots \, \cdot [\bar{\Lambda}^P_{w_s}] = d[\bar{\Lambda}^P_{e}]$ in  $H^*(G/P)$, for some nonzero $d$.
 The converse is not true in general (cf., Theorem ~\ref{T1}).

\begin{definition} \label{basicbundle} Let $w\in W^P$.  Since $T_w:=T_e(\Lambda_w)$ is a $B_L$-module (by (\ref{eqn2})), we have the  $P$-equivariant vector bundle $\mathcal{T}_w:=P \tset{B_L}
T_w$ on $P/B_L$ associated to the principal $B_L$-bundle $P\to P/B_L$
via the $B_L$-module $T_w$. In particular, we have the $P$-equivariant
 vector bundle $\mathcal{T}:=P \tset{B_L} T$ (where, as in Section 2, $
T:=  T_e (G/P)$) and  $\mathcal{T}_w$ is canonically a $P$-equivariant subbundle of $\mathcal{T}$.
Take the top exterior powers $\det (\mt/\mathcal{T}_w)$
and $\det (\mathcal{T}_w)$, which are  $P$-equivariant
 line bundles on $P/B_L$. Observe that, since $T$ is a $P$-module, the
$P$-equivariant  vector bundle $\mathcal{T}$ is $P$-equivariantly isomorphic with the product bundle  $P/B_L \times T$ under the map $\xi: P/B_L \times T
\to \mathcal{T}$ taking $(pB_L,v)\mapsto (p, p^{-1}v)$ mod $B_L$, for
$p\in P$ and $v\in T$; where $P$ acts on $P/B_L \times T$ diagonally. We will
often identify  $ \mathcal{T}$  with the product bundle  $P/B_L \times T$ under  $\xi$.

Similarly, for any $\lambda \in X(H)$, we have a  $P$-equivariant
 line bundle $\cl (\lambda)= \cl_P (\lambda)$  on $P/B_L$ associated to the principal $B_L$-bundle $P\to P/B_L$
via the one dimensional $B_L$-module $\lambda^{-1}$. (As observed in Section
2, any  $\lambda \in X(H)$ extends uniquely to a character of $B_L$.) The twist
in the definition of  $\cl (\lambda)$ is introduced so that the dominant characters correspond to the dominant line bundles.

For $w\in W^P$, define the character $\chi_w\in\frh^*$ by
$$\chi_w=\sum_{\beta\in (R^+\setminus R^+_\fl)\cap w^{-1}R^{+}} \beta \,.$$
 Then, from [Ku2, 1.3.22.3] and (\ref{eqn1}),
\begin{equation}\label{eqn5}
\chi_w = \rho -2\rho^L + w^{-1}\rho ,
\end{equation}
where $\rho$ (resp. $\rho^L$) is half the sum of roots in $R^+$ (resp. in
$ R^+_\fl$).
\end{definition}

The following lemma is easy to establish.

\begin{lemma}\label{mt}
For $w\in W^P$,  as $P$-equivariant
 line bundles  on $P/B_L$, we have:
$$\det (\mt/\mathcal{T}_w)=\cl (\chi_w).$$
Observe that, since $\chi_1$ is a $P$-module (as $T$ is a $P$-module), the line bundle
$\cl (\chi_1)$ is a trivial line bundle. However, as an $H$-equivariant line bundle, it is nontrivial in general as the character $\chi_1$ restricted to the connected center of $L$ is nontrivial.
\end{lemma}

 Let $\mt_s$ be the $P$-equivariant product  bundle
$(P/B_L)^s \times T \to (P/B_L)^s$ under the diagonal action of $P$ on
$(P/B_L)^s \times T $. Then, $\mt_s$ is canonically  $P$-equivariantly
 isomorphic with the pull-back bundle $\pi_j^*(\mathcal{T})$, for any
$1\leq j\leq s$, where   $\pi_j: (P/B_L)^s  \to P/B_L$ is the projection onto the $j$-th factor.  For any  $w_1, \cdots, w_s \in W^P$,
 we have a  $P$-equivariant map of vector bundles on $(P/B_L)^s$:
 \begin{equation}\label{map}
\Theta=\Theta_{(w_1, \cdots, w_s)
}:\mt_s \to \oplus_{j=1}^s
\pi_j^*({\mt}/\mathcal{T}_{w_j})
\end{equation}
obtained as the direct sum of the canonical projections $\mt_s \to
\pi_j^*({\mt}/\mathcal{T}_{w_j})$ under the identification $\mt_s\simeq
\pi_j^*(\mathcal{T})$.
Now, assume that  $w_1, \cdots, w_s \in W^P$ satisfies the condition
(\ref{dim0}). In this case,
 we have
the same rank on the two sides of the map (\ref{map}). Let $\theta$ be the
 bundle map obtained from $\Theta$ by taking the top exterior power:
 \begin{equation}\label{determinant}
\theta=\det(\Theta): \det \mt_s \to \det \Bigl(
\mt/\mathcal{T}_{w_1}\Bigr) \boxtimes \cdots \boxtimes
\det\bigl(\mt/\mathcal{T}_{w_s}\bigr),
  \end{equation}
where $ \boxtimes$ denotes the external tensor product. Clearly, $\theta$
is  $P$-equivariant
 and hence one can view $\theta$ as a  $P$-invariant
element in $$H^0\Biggl( (P/B_L)^{s}, \det (\mt_s )^*\otimes
\Bigl(\det \bigl(\mt/\mt_{w_1})\boxtimes \cdots \boxtimes \det \bigl(\mt/\mt_{w_s}\bigr) \Bigr)\Biggr).$$
\begin{equation}\label{GIT}
=H^0\left( (P/B_L)^{s}, \cl (\chi_{w_1}-\chi_1) \boxtimes \cdots \boxtimes \cl (\chi_{w_s})\right),
\end{equation}
where the above equality follows from Lemma \ref{mt}.

The following lemma is immediate:
\begin{lemma}\label{sectiontheta}
Let $\bar{p}=(\bar{p}_1, \cdots, \bar{p}_s)\in (P/B_L)^s$.
Then, under the assumption (\ref{dim0}),  the following are equivalent:
\begin{enumerate}
\item The restriction of the map $\Theta$ to   the fiber over $\bar{p}$
 is an isomorphism.
\item The section $\theta$ does not vanish at $\bar{p}$.
\item The locally-closed subvarieties $p_1\Lambda_{w_1}^P, \cdots, p_s\Lambda_{w_s}^P$
meet transversally at  $e\in G/P$.
\end{enumerate}
\end{lemma}

The following corollary follows immediately from Proposition \ref{First} and Lemma \ref{sectiontheta}.
\begin{corollary}\label{movL} Let $(w_1, \cdots, w_s)$ be an $s$-tuple of elements of $W^P$ satisfying the condition (\ref{dim0}). Then, we have the following:
\begin{enumerate}
\item The section $\theta$ is nonzero if and only if
$$[\bar{\Lambda}_{w_1}^P]\cdot \,\cdots \,\cdot[\bar{\Lambda}_{w_s}^P]\neq 0\in H^*(G/P).$$
\item The $s$-tuple $(w_1, \cdots, w_s)$ is $L$-movable if and only if the section $\theta$ restricted to $(L/B_L)^s$ is not identically $0$.
 \end{enumerate}
\end{corollary}

\section{Geometric Invariant Theory Revisited}\label{GITT}
We need to consider the Geometric Invariant Theory (GIT) in a nontraditional setting, where a {\it nonreductive} group acts on a {\it nonprojective} variety. To handle
such a  situation, we need to introduce the notion of $P$-admissible one parameter subgroups (cf. Definition \ref{Padmissible}). But first we recall the following definition due to Mumford.
\begin{definition}\label{git} Let $S$ be any (not necessarily reductive) algebraic group
acting on a  (not necessarily projective) variety  $\exx$ and let  $\elal$ be
an $S$-equivariant line bundle on $\exx$.
 Take any $x\in \exx$ and a one parameter subgroup (for short OPS)
 $\lambda \in O(S)$ such that the limit
  $$\lim_{t\to 0}\lambda(t)x$$
exists in $\exx$ (i.e., the morphism ${\lambda}_x:\Bbb{G}_m\to X$ given by
$t\mapsto \lambda(t)x$ extends to a morphism $\tilde{\lambda}_x : \Bbb{A}^1\to X$). This condition is  satisfied by every $x\in \exx$, if $\exx$ is projective.
Then, following Mumford, define a number $\mu^{\elal}(x,\lambda)$ as follows:
 Let $\elal':=\tilde{\lambda}_x^*(\elal)$ be the pull-back line bundle on
 $\Bbb{A}^1$.
  Let $\sigma_1$ be a nonzero vector in the fiber of $\elal$ over $x$. Then, using the $\Bbb{G}_m$-action, $\sigma_1$ extends to
a $\Bbb{G}_m$-invariant section $\sigma$ of $\elal'$ over   $\Bbb{G}_m$. Since $\elal'$ is a line bundle on  $\Bbb{A}^1$ we can speak of the order of vanishing $\mu^{\elal}(x,\lambda)$ of $\sigma$ at $0$ (this number is negative if $\sigma$ has a pole at $0\in \Bbb{A}^1$).

Let $V$ be a finite dimensional representation of $S$ and let
$i:\xx\hookrightarrow \Bbb{P}(V)$ be an $S$-equivariant embedding.
Take $\elal:=i^*(\mathcal O(1))$. Let $\lambda \in O(S)$  and let
$\{e_1,\dots,e_n\}$ be a basis of $V$ consisting of eigenvectors,
i.e.,
 $\lambda(t)\cdot e_l=t^{{\lambda}_l}e_l$, for $l=1,\dots,n$. For any $x\in \xx$ such that  $\lim_{t\to 0}\lambda(t)x$
exists in $\exx$, write $i(x)=[\sum_{l=1}^n x_le_l]$. Then, it is easy to see that, we have ([MFK, Proposition 2.2.3, page 51])
\begin{equation}\label{muvalue}\mu^{\elal}(x,\lambda)=\max_{l: x_l\neq 0}(-{\lambda}_l).
\end{equation}

\end{definition}

We record the following simple properties of $\mu^{\elal}(x,\lambda)$:
\begin{proposition}\label{propn14} For any $x\in \exx$ and $\lambda \in O(S)$ such that $\lim_{t\to 0}\lambda(t)x$
exists in $\exx$, we have the following (for any $S$-equivariant line bundles
$\elal, \elal_1, \elal_2$):
\begin{enumerate}
\item[(a)]
$\mu^{\elal_1\otimes\elal_2}(x,\lambda)=\mu^{\elal_1}(x,\lambda)+\mu^{\elal_2}(x,\lambda).$
\item[(b)] If there exists $\sigma\in H^0(\exx,\elal)^S$ such that $\sigma(x) \neq 0$, then  $\mu^{\elal}(x,\lambda)\geq 0.$
\item[(c)] If $\mu^{\elal}(x,\lambda)=0$, then any element of $H^0(\exx,\elal)^S$ which does not vanish at $x$ does not vanish at $\lim_{t\to 0}\lambda(t)x$ as well.
\item[(d)] For any $S$-variety $\exx'$ together with an $S$-equivariant morphism $f:\exx'\to \exx$ and any $x'\in \exx'$ such that  $\lim_{t\to 0}\lambda(t)x'$
exists in $\exx'$, we have
$$\mu^{f^*\elal}(x',\lambda)=\mu^{\elal}(f(x'),\lambda).$$
\end{enumerate}
\end{proposition}
\begin{definition}\label{Padmissible} We call an OPS $\lambda \in O(P)$,
 {\it $P$-admissible}
if the limit
$\lim_{t\to 0}\lambda(t)x$
exists in $P/B_L$ for every $x\in P/B_L$. When the reference to $P$ is clear from the context, we will abbreviate $P$-admissible by {\it admissible}.

From the conjugacy of $L$ and $B_L$, it is easy to see that the notion of the
admissibility of $\lambda$ does not depend upon the choice of the Levi subgroup $L$ or the Borel subgroup $B_L$ of $L$.
\end{definition}

For an OPS  $\lambda \in O(G)$, set
$$\dot{\lambda}:=\frac{d\lambda(t)}{dt}|_{t=1}\in \frg,$$
its tangent vector.

Also define the {\it associated parabolic subgroup} $P(\lambda)$ of $G$ by
$$P(\lambda):= \{g\in G: \lim_{t\to 0}\lambda(t)g\lambda(t)^{-1} \,
\text{exists in}\,G\}.$$
We also denote $P(\lambda)$ sometimes by $P(\dot\lambda)$. This should not create any confusion since $\dot\lambda$ uniquely determines $\lambda$.

The following lemma gives a characterization of admissible one parameter subgroups in $P$.

\begin{lemma}\label{lemmar} Let $\lambda \in O(P)$. Then, $\lambda$ is admissible iff
$$P(\lambda) \supset U.$$
Equivalently,  $\lambda$ is admissible iff there exists $p\in P$ such that
$\lambda_o:=p\lambda p^{-1}$ lies in $H$ and $\beta(\dot{\lambda_o}) \geq 0$ for all $\beta \in R^+\setminus R^+_\mathfrak l.$
\end{lemma}
\begin{proof}
Clearly, $\lambda$ is admissible iff any conjugate $p\lambda p^{-1}$
is admissible, for $p\in P$. Moreover, $P(p\lambda p^{-1})=pP(\lambda ) p^{-1}.$ Thus, $U$ being normal in $P$, the condition that
$P(\lambda) \supset U$  is invariant under the conjugation of $\lambda$ via $p$.
Further, since any OPS in $P$ can be conjugated via an element $p\in P$
inside the maximal torus $H$, we can assume that $\Imo \lambda \subset H$.
Decompose
$$P/B_L \simeq  U \times L/B_L\,\,\text{as varieties}.$$
For any OPS $\lambda$ in $H$ (in particular in $L$), and any
$u\in U, l\in L,$
$$\lim_{t\to 0}\lambda(t)(ulB_L)=\lim_{t\to 0}\bigl(\lambda(t)u\lambda(t)^{-1}
(\lambda(t)lB_L)\bigr).$$
But, since $L/B_L$ is projective, the limit on the left side exists iff
$\lim_{t\to 0}\lambda(t)u\lambda(t)^{-1}$ exists in $U$, which is equivalent to the condition that $P(\lambda) \supset U.$  This proves the first part of the lemma.

The second part follows readily from the first part.

\end{proof}

\begin{definition}\label{mumfordindex}  Let $x=u lB_L\in P/B_L$, for $u\in U$ and $l\in L$. Let
 $\lambda$ be an admissible OPS
 and let $P_L(\lambda):= P(\lambda)\cap L.$ Then,  $P_L(\lambda)$ is a parabolic subgroup of $L$. To see this, write  $\lambda = u'\lambda_o {u'}^{-1}$,
for $u'\in U$ and with $\lambda_o$ an OPS in $L$. This gives $P(\lambda)= u'P(\lambda_o){u'}^{-1}$. But since $\lambda$ is admissible by assumption, and hence so is $\lambda_o$; and thus by Lemma  \ref{lemmar},
 $P(\lambda_o)\supset U$. This gives   $P(\lambda) = P(\lambda_o)$. But, as is well known, $P_L(\lambda_o)$ is a parabolic subgroup of $L$ and hence so is
$P_L(\lambda)$.

Write $P_L(\lambda)=l_1Ql_1^{-1}$, for some $l_1\in L$,
 where $Q$ is a standard parabolic subgroup of $L$, i.e., $Q\supset B_L$.
Let $l^{-1}l_1 \in B_LwQ$, where $w\in W_L/W_Q$. Then, clearly $w\in
W_L/W_Q$ does not depend upon the choices of the representatives $l$ (in
$lB_L$) and $l_1$. We define the {\it relative position}
$[x,\lambda]$  to be $w \in W_L/W_Q$. This  satisfies:
\begin{equation}\label{eqn10}
[px,p\lambda p^{-1}]=[x,\lambda],\,\text{ for any }\,p\in P.
\end{equation}

Observe that, if $\lambda$ is an OPS lying in the center of $L$, then for any
$x\in P/B_L$,
\begin{equation}\label{eqn11}
[x,\lambda]=1.
\end{equation}
 For an OPS $\lambda$ in  $P$, we can choose $p\in P$ such that $\lambda_o:=
p\lambda p^{-1}$ is an OPS in $H$ and, moreover, $\dot{\lambda}_o\in \mathfrak h$ is
$L$-dominant (i.e., $\alpha_i(\dot{\lambda}_o)\geq 0$ for all $\alpha_i\in
\Delta(P)$). Set,
\begin{equation}\label{Xtilde}
{X}_{\lambda}^P= \dot{\lambda}_o.
\end{equation}
Then, ${X}_{\lambda}^P$ is well defined. Moreover, $W_Q$ fixes $
{X}_{\lambda}^P$. We shall abbreviate ${X}_{\lambda}^P$ by
${X}_{\lambda}$ when the reference to $P$ is clear.

The following lemma is a generalization of the corresponding result in [BeSj, Section 4.2].
\end{definition}
\begin{lemma}\label{l1}Let $\lambda$ be an admissible $OPS$ in $P$, $x
\in P/B_L$ and
$\chi\in X(H)$. Then, we have the  following formula:
$$\mu^{\cl(\chi)}(x,\lambda)=-\chi([x,\lambda] {X}_{\lambda}).$$
\end{lemma}
\begin{proof}
Let $p\in P$ be such that $\lambda= p\lambda_op^{-1}$ (where $\lambda_o$ is an OPS in $H$ such that $\dot{\lambda}_o=X_\lambda$). Let $\bar{x} \in P$ be a lift of $x$.  We seek an OPS
$b(t)$ in $B_L$  so that $\lambda_o(t)p^{-1}\bar{x}b(t)$ has a limit in $P$ as  $t\to 0$. Let
$p=u_1l_1,\bar{x}=ul$, for $u,u_1\in U$ and $l,l_1\in L$,  and let
$Q$ be the standard parabolic subgroup of $L$ as in Definition ~\ref{mumfordindex}. Choose $w \in W_L/W_Q, b_l \in [B_L,B_L]$ and $q\in Q$ so that
 $l_1=lb_lwq$. Then, by the definition,
\begin{equation}\label{eqn10}
w=[x,\lambda].
\end{equation}
Thus, for any $b(t) \in B_L$,
$$\lambda_o(t)p^{-1}\bar{x}b(t)=\lambda_o(t)l_1^{-1}u_1^{-1}u l_1 l_1^{-1}lb(t)
=\lambda_o(t)\hat{u}\lambda_o(t)^{-1} \lambda_o(t)q^{-1}\lambda_o(t)^{-1}\lambda_o(t)w^{-1}b_l^{-1}b(t),$$
where $\hat{u}:=l_1^{-1}u_1^{-1}u l_1\in U$.
In view of Lemma~\ref{lemmar} and since $P_L(\lambda_o)=Q$,   we find that the OPS
$$b(t):=b_lw\lambda_o(t)^{-1}w^{-1}$$
``works'', i.e.,   $\lambda_o(t)p^{-1}\bar{x}b(t)$ has a limit in $P$ as
$t\to 0$.

Consider the  $\Bbb{G}_m$-invariant section $\sigma(t):=(\lambda(t)\bar{x},1)$\,mod $B_L$ of $\lambda_x^*(\cl_\chi)$ over  $\Bbb{G}_m$, where, as in
Definition~\ref{git}, $\lambda_{{x}}: \Bbb{G}_m \to P/B_L$ is the map $t\mapsto \lambda(t)\cdot {x}$. Then, the section $\sigma(t)$ corresponds to the function $\Bbb{G}_m \to \Bbb A^1, \,t\mapsto \chi^{-1}(b(t)^{-1})$. From this we get the lemma by using (\ref{eqn10}).
\end{proof}

\section{Criterion for $L$-movability}\label{ref11}
The aim of this section is to prove the following characterization of  $L$-movability.
\begin{theorem}\label{T1}
Assume that  $(w_1, \cdots, w_s)\in (W^P)^s$  satisfies equation (\ref{dim0}).
Then, the following are equivalent.

(a) $(w_1, \cdots, w_s)$ is $L$-movable.

(b) $[\bar{\Lambda}^P_{w_1}]\cdot\, \cdots \,\cdot [\bar{\Lambda}^P_{w_s}] = d[\bar{\Lambda}^P_{e}]$ in  $H^*(G/P)$, for some nonzero $d$,
 and for each  $\alpha_i\in \Delta\setminus\Delta(P)$, we have
$$\bigl((\sum_{j=1}^s\,\chi_{w_j})-\chi_1\bigr)(x_{i})=0,$$
where $\chi_w$ is as defined in Definition \ref{basicbundle} and $x_i\in \fh$ is defined by (\ref{eqn0}).
\end{theorem}

\begin{proof} (a)$\Rightarrow$(b):
Let $(w_1,\cdots, w_s)\in (W^P)^s$ be $L$-movable. Consider the restriction
$\hat{\theta}$   of the $P$-invariant section
$$\theta\in H^0\left( (P/B_L)^{s}, \cl (\chi_{w_1}-\chi_1) \boxtimes
\cdots \boxtimes \cl (\chi_{w_s})\right)$$
to $(L/B_L)^s$, where $\theta$ is as defined in equations
(\ref{determinant})--(\ref{GIT}). Then, $\hat{\theta}$  is non-vanishing by Corollary \ref{movL}. But for
$$H^0\left( (L/B_L)^{s},  \cl (\chi_{w_1}-\chi_1) \boxtimes
\cdots \boxtimes \cl (\chi_{w_s})\right)^L$$
 to be non-empty, the center of $L$ should act trivially (under the diagonal action)
on  $\cl (\chi_{w_1}-\chi_1) \boxtimes
\cdots \boxtimes \cl (\chi_{w_s})$
restricted to $(L/B_L)^s$. This gives
$$\sum_{j=1}^s\chi_{w_j}(h)=\chi_1(h),$$
for all $h$ in the Lie algebra $\mathfrak z_L$ of the center of $L$.
For any $\alpha_i \in \Delta \setminus\Delta(P)$, $h=x_i$ clearly lies in
 $\mathfrak z_L$. Taking $h=x_{i}$ in the above equation, we obtain the implication
(a)$\Rightarrow$(b)  by using Corollary \ref{movL}.

We now prove the implication (b)$\Rightarrow$(a). Thus, we are given that
 for each $\alpha_i\in \Delta\setminus \Delta(P)$,
$$\bigl((\sum_{j=1}^s\chi_{w_j})-\chi_1\bigr)(x_i)=0.$$
Moreover, by Corollary \ref{movL},  $\theta(\bar{p}_1, \cdots, \bar{p}_s)\neq 0$,
for some $\bar{p}_j\in P/B_L$. Consider the central OPS of $L$:
$\lambda (t):=\prod_{\alpha_i\in\Delta\setminus\Delta(P)} t^{x_i}.$
Clearly,
\begin{equation}\label{eqn12}
\dot{\lambda}=\sum_{\alpha_i\in\Delta\setminus\Delta(P)}x_i.
\end{equation}
Thus, by Lemma \ref{lemmar},
$\lambda$ is admissible.
 For any $x=ulB_L\in P/B_L$, with $u\in U$ and $l\in L$,
$$\lim_{t\to 0}\lambda(t)x=\lim_{t\to 0}\lambda(t)u\lambda(t)^{-1}
(\lambda(t)l) B_L.$$
But, since $\beta (\dot\lambda) >0$, for all $\beta \in R^+\setminus R^+_\mathfrak l$, we get
$$\lim_{t\to 0}\lambda(t)u\lambda(t)^{-1}=1.$$
Moreover, since  $\lambda(t)$ is central in $L$,
  the limit
$\lim_{t\to 0}\lambda(t)lB_L$ exists and equals $lB_L$. Thus,
$\lim_{t\to 0}\lambda(t)x$ exists and lies in $L/B_L$.

Now, let $\elal$ be the $P$-equivariant line bundle $\cl (\chi_{w_1}-\chi_1) \boxtimes
\cdots \boxtimes \cl (\chi_{w_s})$  on $\exx:=(P/B_L)^s$, and $\bar{p}:=
(\bar{p}_1,\cdots,\bar{p}_s)\in \exx$. Then, by Lemma
\ref{l1} and equations (\ref{eqn11}) and (\ref{eqn12}), we get
$$\begin{aligned}
\mu^{\elal}(\bar{p},\lambda)&=(\chi_{1}-\chi_{w_1})
([\bar{p}_1, \lambda] \dot{\lambda})- \sum_{j=2}^s\chi_{w_j}
([\bar{p}_j, \lambda] \dot{\lambda})\\
&=-\sum_{\alpha_i\in\Delta\setminus\Delta(P)}\bigl(
\bigl(\bigl(\sum_{j=1}^s\chi_{w_j}\bigr)-\chi_1\bigr)(x_i)\bigr)\\
&=0,\,\,\text{by assumption}.
\end{aligned}$$
Therefore, using Proposition \ref{propn14}(c) for $S=P$,   $\theta$ does not vanish at
$\lim_{t\to 0}\lambda(t)\bar{p}.$
But, from the above, this limit exists as an element of $(L/B_L)^s$.
Hence, $(w_1, \cdots, w_s)$ is $L$-movable by Corollary ~\ref{movL}.

\end{proof}

\section{Deformation of Cup Product in $H^*(G/P)$}

 Define the structure constants $c^{w}_{u,v}$ under the intersection product
in $H^*(G/P, \Bbb Z)$ via  the formula
\begin{equation}
[\bar{\Lambda}^P_u]\cdot[\bar{\Lambda}^P_v]=\sum_{w\in W^P} c^w_{u,v}[\bar{\Lambda}^P_w].
\end{equation}
The number $c^w_{u,v}$ is the number of points (counted with multiplicity) in the intersection  $g_1[\bar{\Lambda}^P_u]\cap g_2[\bar{\Lambda}^P_v]\cap g_3[\bar{\Lambda}^P_{w_oww^P_o}]$
for generic $(g_1,g_2,g_3)\in G^3$. If the generic intersection is infinite, we set $c^w_{u,v}=0$.

Introduce the indeterminates ${\defpar}_i$ for each $\alpha_i\in \Delta\setminus\Delta(P)$ and write a deformed cup product
\begin{equation}\label{EE}
[\bar{\Lambda}^P_u]\odot[\bar{\Lambda}^P_v]=\sum_{w\in W^P} \bigl(\prod_{\alpha_i\in \Delta\setminus\Delta(P)}
{\defpar}_i^{(\chi_w-(\chi_u+\chi_v))(x_i)} \bigr)c^w_{u,v}[\bar{\Lambda}^P_w],
\end{equation}
where $\chi_w$ is defined in Definition \ref{basicbundle}. Extend this to a $\Bbb Z[{\defpar}_i]_{\alpha_i\in \Delta\setminus\Delta(P)}-$linear product structure on $H^*(G/P, \Bbb Z) \otimes_\Bbb Z \Bbb Z[{\defpar}_i].$ By the next Proposition \ref{product}, if $c^w_{u,v}\neq 0$,
 $$(\chi_w-(\chi_u+\chi_v))(x_i)\geq 0, \,\,\,\text{ for any}\,
\alpha_i\in \Delta\setminus\Delta(P)$$
 and hence the above product indeed lies in $H^*(G/P, \Bbb Z) \otimes_\Bbb Z \Bbb Z[{\defpar}_i].$
Clearly, $\odot$ is commutative. {\it This product should not be confused with the 
product in the quantum cohomology of $G/P$.}

Recall the definition of $T, T_w$ from Section 2.

\begin{lemma}\label{long}
For any $w\in W^P,$

(a)
 $T_w\oplus w^P_o (T_{w_oww^P_o})=T$.

(b)  $\chi_w+w^P_o\chi_{w_oww^P_o}=\chi_1$.

(c) $(\chi_w+\chi_{w_oww^P_o}-\chi_1)(x_i)=0$, for all $\alpha_i\in \Delta\setminus\Delta(P)$.

(d) For any $v,w\in W^P$ such that $\ell (v)=\ell (w)$,  $[\bar{\Lambda}^P_v]\odot[\bar{\Lambda}^P_{w_oww^P_o}]=\delta_{v,w} [\bar{\Lambda}^P_{1}].$
\end{lemma}
\begin{proof}
As in Section 2, let $R(T_w)$ denote the set of roots $\alpha$ such that the root space $\mathfrak g_\alpha \subset T_w$, i.e.,
$T_w=\oplus_{\alpha\in R(T_w)}\,\mathfrak g_\alpha$. Then,
\begin{equation}\label{eqn15}
 R(T_w) = w^{-1}R^+\cap \bigl(R^-\setminus R^-_{\mathfrak l}\bigr).
\end{equation}
To prove (a), it suffices to show that we have the disjoint union:
\begin{equation}\label{eqn16}
 R^-\setminus R^-_{\mathfrak l}=R(T_w) \sqcup w_o^P\cdot  R(T_{w_oww^P_o}).
\end{equation}
 Since  $w_o^P$ keeps $ R^-\setminus R^-_{\mathfrak l}$ stable,
by (\ref{eqn15}),
$$\begin{aligned}
w_o^P\cdot  R(T_{w_oww^P_o}) &=w^{-1}w_oR^+\cap
\bigl(R^-\setminus R^-_{\mathfrak l}\bigr)\\
&=w^{-1}R^-\cap \bigl(R^-\setminus R^-_{\mathfrak l}\bigr).
\end{aligned}$$
From this and equation (\ref{eqn15}), equation (\ref{eqn16}) follows.

To prove (b), by equation (\ref{eqn5}),
$$\begin{aligned}
\chi_w+w_o^P\chi_{w_oww_o^P}&=\rho-2\rho^L+w^{-1}\rho+w_o^P\rho-2w_o^P\rho^L+w^{-1}w_o\rho\\
&=\rho-2\rho^L+w^{-1}\rho+w_o^P(\rho-\rho^L)-w_o^P\rho^L-w^{-1}\rho\\
&=\rho-2\rho^L+\rho-\rho^L+\rho^L,\,\,\text{since}\, w_o^P\,\,\text{permutes}\,
R^+\setminus R^+_{\mathfrak l}\\
&=2\rho - 2\rho^L\\
&=\chi_1.
\end{aligned}$$
This proves (b).

For any $\alpha_i \in \Delta \setminus \Delta(P)$,  $x_i$ is central in $\mathfrak l$; in particular, $w_o^P$ acts trivially on $x_i$. Thus (c) follows from (b).

By [KLM, Lemma 2.9],
\begin{equation}\label{eqn17}
[\bar{\Lambda}^P_v]\cdot[\bar{\Lambda}^P_{w_oww^P_o}]=\delta_{v,w} [\bar{\Lambda}^P_{1}].
\end{equation}
Thus, (d) follows from the defining equation (\ref{EE}) and the (c)-part.

\end{proof}
\begin{proposition}\label{product}

(a) For any $u,v,w\in W^P$ such that $c^w_{u,v}\neq 0$, we have
\begin{equation}\label{eqn18}
 (\chi_w-\chi_u-\chi_v)(x_i)\geq 0,\,\,\text{ for each}\,\, \alpha_i\in \Delta\setminus\Delta(P).
\end{equation}

(b) The product $\odot$  in $H^*(G/P)\otimes\Bbb{Z}[{\defpar}_i]$  is associative.

(c) For $(w_1,\dots,w_s)\in (W^P)^s$ and $w\in W^P$, the coefficient of $
[\bar{\Lambda}^P_{w}]$ in
$$[\bar{\Lambda}^P_{w_1}]\odot\dots\odot[\bar{\Lambda}^P_{w_s}]$$
is $$\prod_{\alpha_i\in \Delta\setminus\Delta(P)}
{\defpar}_i^{(\chi_w-\sum_{j=1}^s\chi_{w_j})(x_i)}$$
times the coefficient of 
$[\bar{\Lambda}^P_{w}]$ in the usual cohomology product
$$[\bar{\Lambda}^P_{w_1}]\cdot\dots\cdot[\bar{\Lambda}^P_{w_s}].$$
  \end{proposition}
  \begin{proof}
 By equation (\ref{eqn17}), $c^w_{u,v}\neq 0$ iff
$$[\bar{\Lambda}^P_u]\cdot[\bar{\Lambda}^P_v]\cdot[\bar{\Lambda}^P_{w_oww_o^P}]=d[\bar{\Lambda}^P_1],
$$ for some nonzero $d$. Thus, by Theorem \ref{T2}(a), for any
$\alpha_i \in \Delta \setminus \Delta(P)$,
$$(\chi_u+\chi_v+\chi_{w_oww_o^P}-\chi_1) (x_i)\leq 0.$$
By Lemma \ref{long}(c), this gives
$$(\chi_u+\chi_v-\chi_{w}) (x_i)\leq 0,$$
proving (a).

To prove (b), write
$$\begin{aligned}
\bigl([\bar{\Lambda}^P_u]\odot [\bar{\Lambda}^P_v]\bigr) \odot [\bar{\Lambda}^P_w]&=\sum_{\theta\in W^P} \bigl(\prod
{\defpar}_i^{(\chi_\theta-(\chi_u+\chi_v))(x_i)} \bigr)c^\theta_{u,v}[\bar{\Lambda}^P_\theta] \odot [\bar{\Lambda}^P_w]\\
&=\sum_{\theta,\eta} \bigl(\prod
{\defpar}_i^{(\chi_\eta-\chi_w-\chi_u-\chi_v)(x_i)} \bigr)c^\theta_{u,v} c^\eta_{\theta,w}[\bar{\Lambda}^P_\eta].
\end{aligned}$$
Similarly,
$$
[\bar{\Lambda}^P_u]\odot \bigl([\bar{\Lambda}^P_v] \odot [\bar{\Lambda}^P_w]\bigr)
=\sum_{\theta,\eta} \bigl(\prod
{\defpar}_i^{(\chi_\eta-\chi_u-\chi_v-\chi_w)(x_i)} \bigr)c^\eta_{u,\theta} c^\theta_{v,w}[\bar{\Lambda}^P_\eta].
$$
So, the associativity of $\odot$ follows from the corresponding associativity in $H^*(G/P)$ under the standard cup product. The property (c) is immediate from the definitions and property (b).
  \end{proof}

\begin{definition} \label{levicohomology} The  cohomology of $G/P$
 obtained by setting each ${\defpar}_i=0$ in
$(H^*(G/P, \Bbb Z)\otimes\Bbb{Z}[{\defpar}_i],\odot)$ is denoted by
$(H^*(G/P, \Bbb Z),\odot_0)$. Thus, as a $\Bbb Z$-module, this  is the same as the singular cohomology
 $H^*(G/P, \Bbb Z)$. This has essentially the effect of ignoring
all the non $L$-movable intersections. By the above proposition,  $(H^*(G/P, \Bbb Z),\odot_0)$ is associative (and commutative).
Moreover, by Lemma \ref{long}(d), it continues to satisfy the Poincar\'e duality.
\end{definition}
\begin{lemma} \label{minuscule} Let $P$ be a minuscule maximal standard parabolic subgroup of $G$ (i.e., the simple root $\alpha_p \in \Delta \setminus \Delta(P)$ appears with coefficient
$1$ in the highest root of $R^+$ and hence in all the roots in $R(\fu_P)$). Then, for any $u,v \in W^P$,
$$[\bar{\Lambda}^P_u]\odot[\bar{\Lambda}^P_v]=[\bar{\Lambda}^P_u]\cdot[\bar{\Lambda}^P_v].$$
\end{lemma}
\begin{proof} By the definition of $\odot$, it suffices to show that for any $w\in W^P$ such that $c_{u,v}^w\neq 0$,
\begin{equation} \label{eqn5.1}
 (\chi_w-(\chi_u+\chi_v))(x_p)= 0.
\end{equation}
By the definition of $\chi_w$ (cf. Definition \ref{basicbundle}), since $P$ is minuscule,
\begin{equation} \label{eqn5.2}
 \chi_w (x_p)=  \mid w^{-1}R^+\cap \bigl(R^+\setminus R^+_{\mathfrak l}\bigr)\mid
=\codim(\Lambda^{P}_w;G/P),
\end{equation}
where the last equality follows from equation (\ref{eqn15}). Moreover, since
 $c_{u,v}^w\neq 0$,
\begin{equation} \label{eqn5.3}
\codim(\Lambda^{P}_u;G/P)+\codim(\Lambda^{P}_v;G/P)
=\codim(\Lambda^{P}_w;G/P).
\end{equation}
Combining equations (\ref{eqn5.2}) and (\ref{eqn5.3}), we get equation (\ref{eqn5.1}).
\end{proof}

\section{Solution of the Eigenvalue Problem}
Our aim in this section is to prove that for the solution of the eigenvalue problem one may restrict to those inequalities coming from $L$-movable intersections with intersection number one.
\subsection{Principal criterion for the nontriviality of the space of invariants}\label{GIT1}
Let $\nu_1, \cdots, \nu_s\in X(H)$ be dominant weights and let $\elal$ be the $G$-linearized line bundle $\cl(\nu_1)\boxtimes \cdots \boxtimes\cl(\nu_s)$ on $(G/B)^s$.

In this subsection, we give our  principal criterion to decide whether there exists an integer $N>0$ such that
$H^0((G/B)^s,\Bbb{L}^N)^G \neq 0.$

Let $P_1, \cdots, P_s$ be the standard parabolic subgroups such that the line bundle
$\cl(\nu_j)$ on $G/B$ descends as an ample line bundle on $G/P_j$, still
denoted by $\cl(\nu_j)$. Consider $\xx :=G/P_1 \times \cdots \times G/P_s$ with the diagonal action of $G$ and  $\elal$  the $G$-linearized ample line bundle $\cl(\nu_1)\boxtimes \cdots \boxtimes\cl(\nu_s)$ on $\xx$. Let $\pi: (G/B)^s \to \xx$  be the canonical projection map.
\begin{lemma}\label{phone}
 Let $x\in (G/B)^s$ and let $x'=\pi(x)$. The following are equivalent.
\begin{enumerate}
\item[(a)]  For some integer $N>0$, there exists
$\sigma\in H^0((G/B)^s,\Bbb{L}^N)^G $ such that $\sigma(x)\neq 0$, where
 $G$ acts diagonally.
\item[(b)] For every OPS $\lambda \in O(G)$, $\mu^\elal(x,\lambda)\geq 0$.
\item[(c)] For every OPS $\lambda \in O(G)$, $\mu^{\elal}(x',\lambda)\geq 0$.
\item[(d)]  For some integer $N>0$, there exists $\sigma'\in H^0(\exx,\Bbb{L}^N)^G$ such that $\sigma'(x')\neq 0$, i.e., $x'$ is a semistable point of $\xx$
(with respect to the $G$-linearized ample line bundle $\elal$).
\end{enumerate}
\end{lemma}
\begin{proof}
We first note that for any integer $N>0$, the map $\pi$ induces an
 isomorphism of $G$-modules:
\begin{equation}\label{equalitie}
H^0(\exx,\Bbb{L}^N) \simeq H^0((G/B)^s,\Bbb{L}^N).
\end{equation}
It follows from Proposition ~\ref{propn14}(d) that (b) is equivalent to
(c).

By Hilbert-Mumford theorem [MFK, Theorem 2.1], since $\elal$ is ample on $\exx$, (c) is equivalent to (d).

By equation (\ref{equalitie}), (d) is equivalent to (a) and we are done.
\end{proof}

We now state one of the main theorems of this paper.
\begin{theorem}\label{EVT} Let $G$ be a connected semisimple group.
With the notation as above, the following are equivalent:
\begin{enumerate}
\item[(i)] For some integer $N>0$,
$$H^0((G/B)^s,\Bbb{L}^N)^G \neq 0. $$
\item[(ii)]For every standard maximal parabolic subgroup $P$ in $G$ and every choice of
{\bf $L$-movable}
$s$-tuples  $(w_1, \cdots, w_s)\in (W^P)^s$ such that
$$[\bar{\Lambda}^P_{w_1}]\cdot\, \cdots \,\cdot [\bar{\Lambda}^P_{w_s}] = [\bar{\Lambda}^P_{e}] \in H^*(G/P,\Bbb{Z}),$$
  the following inequality holds:
\begin{equation}\label{eqn29}
\sum_{j=1}^s \nu_j(w_jx_{i_P})\leq 0,
\end{equation}
where $\alpha_{i_P}$ is the simple root in $\Delta\setminus \Delta (P).$
\item[(iii)]For every standard maximal parabolic subgroup $P$ in $G$ and every choice of
$s$-tuples  $(w_1, \cdots, w_s)\in (W^P)^s$ such that
$$[\bar{\Lambda}^P_{w_1}]\odot_0\, \cdots \,\odot_0 [\bar{\Lambda}^P_{w_s}] = [\bar{\Lambda}^P_{e}] \in \bigl(H^*(G/P,\Bbb{Z}), \odot_0\bigr),$$
  the above inequality (\ref{eqn29}) holds.
\end{enumerate}
\end{theorem}
The proof of this theorem will be given in subsection 7.3.
\begin{remark} (a) The above theorem remains true for any connected reductive $G$ 
provided we assume that $\sum_j{\nu_j}_{\vert \mathfrak z(\mathfrak g)}=0$, where 
$\mathfrak z(\mathfrak g)$ is the center of $\mathfrak g$. 

(b) The implication ($i$) $\Rightarrow$ ($ii$) (and hence also the implication ($i$) 
$\Rightarrow$ ($iii$)) in the above theorem rermains true for any (not necessarily maximal) standard parabolic
$P$ of $G$. 
\end{remark}  
\subsection{Maximally destabilizing one parameter subgroups}\label{Kempf}
We recall the definition of the Kempf's OPS attached to an unstable point, which is in some sense `most destabilizing' OPS. The exposition  below follows the paper of Hesselink [H].

Let $X$ be a projective variety with the action of a reductive group $G$ and let
$\mathcal{L}$ be  a $G$-linearized ample line bundle on $X$.

We introduce the set $M(G)$ of fractional OPS in $G$. This is the set
 consisting of the ordered pairs $(\delta,a)$, where $\delta:\gm\to G$ is an OPS of $G$ and $a\in \Bbb{Z}_{>0}$, modulo the equivalence relation $(\delta,a)\simeq(\nu,b)$ if $\delta^b=\nu^a$. An OPS $\delta$ of $G$ gives the element $(\delta,1)\in M(G)$. The Killing form induces a norm $q$ on $M(G)$ satisfying
$aq(\delta,a)=q(\delta,1)$, and $q(\delta,1)=\Vert \dot{\delta}\Vert$.

We can extend the definition of  $\mu^{\mathcal{L}}(x,\lambda)$ for any element $\lambda=(\delta,a)\in M(G)$ and $x\in X$ by setting $\mu^{\mathcal{L}}(x,\lambda)=
\frac{\mu^{\mathcal{L}}(x,\delta)}{a}$.

For any OPS   $\lambda$ of $G$, recall the definition of the  associated  parabolic subgroup $P(\lambda)$ of $G$ from Definition
\ref{Padmissible}. We extend the definition of $P(\lambda)$ for any $\lambda=(\delta,a)\in M(G)$ by setting
$$P(\lambda)=P(\delta).$$
Then, the group
$$L(\lambda):=\{p\in G\mid \lim_{t\to 0} \lambda(t)p\lambda(t)^{-1}=p\}$$
is a Levi subgroup of $P(\lambda)$ and, moreover,
$$U(\lambda):=\{p\in G\mid \lim_{t\to 0} \lambda(t)p\lambda(t)^{-1}=e\in G\}$$
is the unipotent radical of $P(\lambda)$.

We note the following elementary property: If $\lambda\in M(G)$ and $p\in P(\lambda)$
then
\begin{equation}\label{translate}
\mu^{\mathcal{L}}(x,\lambda)=\mu^{\mathcal{L}}(x,p\lambda p^{-1}).
\end{equation}
For any unstable point $x\in X$, define
$$q^*(x)=\inf_{\lambda\in M(G)}\{q(\lambda)\mid \mu^{\mathcal{L}}(x,\lambda)\leq -1\},$$
and the {\it optimal class}
$$\Lambda(x)=\{\lambda\in M(G)\mid \mu^{\mathcal{L}}(x,\lambda)\leq -1, q(\lambda)=q^*(x)\}.$$
Any $\lambda \in \Lambda(x)$ is called {\it Kempf's OPS associated to} $x$.

The following theorem is due to Kempf (cf. [K, Lemma 12.13]).
\begin{theorem}\label{kempf}
For any  unstable point  $x \in X$ and $\lambda_1$, $\lambda_2\in \Lambda(x)$,
$P(\lambda_1)=P(\lambda_2)$. Moreover, there exist $p_1,p_2\in P(\lambda_1)$
 so that $p_1\lambda_1p_1^{-1}=p_2\lambda_2p_2^{-1}$.

Conversely, for  $\lambda\in \Lambda(x)$ and
$p\in P(\lambda),$ we have $ p\lambda p^{-1}\in \Lambda(x)$ by equation (\ref{translate}).
\end{theorem}
The parabolic $P(\lambda)$ for $\lambda\in \Lambda(x)$ will be denoted by $P(x)$ and called the {\it Kempf's parabolic associated to the unstable point} $x$.
We recall the following theorem due to Ramanan-Ramanathan [RR].
\begin{theorem}\label{RR}
For any  unstable point  $x \in X$ and $\lambda=(\delta,a)\in \Lambda(x)$, let
 $$x_o=\lim_{t\to 0}\,\delta(t)\cdot x \in X.$$
Then, $x_o$ is unstable and $\lambda\in\Lambda(x_o)$.
\end{theorem}
\begin{remark}
The above results are valid  for points $x\in (G/B)^s$ (for the
linearization $\elal$ as in subsection 7.1). We may see this by
applying Theorems ~\ref{kempf}, ~\ref{RR} to  the image
$x'\in\exx$  (and using Proposition \ref{propn14}(d)).
\end{remark}

\subsection{Proof of Theorem ~\ref{EVT}}
We return to the notation and assumptions of subsection 7.1. In particular,
$\nu_1, \cdots, \nu_s\in X(H)$ are dominant weights and $\elal$ denotes the $G$-linearized line bundle $\cl(\nu_1)\boxtimes \cdots \boxtimes\cl(\nu_s)$ on
 $(G/B)^s$.
Also,  $P_1, \cdots, P_s$ are the standard parabolic subgroups such that
$\elal$ descends as an ample line bundle (still denoted by) $\elal$ on
$\xx :=G/P_1 \times \cdots \times G/P_s$. We call a point
 $x\in (G/B)^s$ {\it semistable} (with respect to, not necessarily ample,
$\elal$) if its image in $\xx$ under the canonical map $\pi: (G/B)^s
\to \xx$ is semistable.

By Lemma ~\ref{phone} and equation (\ref{equalitie}), condition (i) of Theorem
\ref{EVT} is equivalent to the following condition:
\begin{enumerate}
\item[($iv$)]The set of semistable points of
$(G/B)^s$ with respect to $\Bbb{L}$ is non empty.
\end{enumerate}

\vskip1ex

Moreover, by the definition of $\odot_0$ and Theorem \ref{T1}, the conditions $(ii)$ and
$(iii)$ of Theorem \ref{EVT} are equivalent.

{\bf Proof of the implication $(iv)\Rightarrow (ii)$ of Theorem \ref{EVT}:} Let
$x=(\bar{g}_1,\dots,\bar{g}_s)
\in (G/B)^s$ be a semistable point, where $\bar{g}_j=g_jB$. Since the set of semistable points is clearly open, we can choose a generic enough $x$ such that the intersection $\cap \,g_jBw_jP$ itself is nonempty. (By assumption,
 $\cap\,  \overline{ g_jBw_j P}$ is nonempty for any $g_j$.) Pick  $f\in
\cap \,g_jBw_jP$. Consider the OPS $\lambda=t^{x_{i_P}}$ which is central in $L$. Since
$P(\lambda) \supset L$ and, clearly, $P(\lambda) \supset B$, we have
$P(\lambda) \supset P$. But, by assumption,  $P$ is a maximal parabolic
subgroup and hence $P(\lambda) =P$ (it is easy to see that $P(\lambda) \neq G$). Moreover, from  Definition \ref{mumfordindex} applied to the case $P=G$, it is easy to see that $[\bar{g}_j, f\lambda f^{-1}]=w_j$ and,
$X_{f\lambda f^{-1}}=x_{i_P}$. Thus, applying  Lemma \ref{l1} for $P=G$, the required inequality (\ref{eqn29}) is
the same as $\mu^\elal(x,f \lambda f^{-1})\geq 0$, but this follows from
Lemma \ref{phone}, since $x$ is semistable by assumption.

\vskip1ex

Before we come to the  proof of the implication $(ii)\Rightarrow (i)$ in Theorem \ref{EVT}, we need to recall the following result due to Leeb-Millson.

Suppose that $x=(\bar{g}_1,\dots,\bar{g}_s)
\in (G/B)^s$ is an unstable point and $P(x)$ the Kempf's parabolic associated
to $x$. Let $\lambda=(\delta, a)$ be a Kempf's OPS associated to  $x$.
Express $\delta(t)=f\gamma(t)f^{-1}$, where $\dot{\gamma}\in \frh_{+}$.
Then $P(\gamma)$ is a standard parabolic. Let $P$ be a maximal parabolic containing $P(\gamma)$. Define $w_j\in W/W_{P(\gamma)}$ by
$fP(\gamma)\in g_j Bw_jP(\gamma)$ for $j=1,\dots,s$. We recall the following
 theorem due to
Leeb-Millson   applicable to any unstable point $x\in (G/B)^s$. We postpone
the proof of this theorem to the next subsection.
\begin{theorem}\label{jm}
\begin{enumerate}
\item[(a)]  The intersection $\bigcap_{j=1}^s  g_jBw_j P\subset G/P $ is the singleton $\{fP\}$.
\item[(b)] For the simple root   $\alpha_{i_P}\in \Delta\setminus \Delta(P)$,
$\sum_{j=1}^s \nu_j(w_j x_{i})>0.$
\end{enumerate}
\end{theorem}

 Now, we come to the  proof of the implication $(ii)\Rightarrow (i)$ in Theorem \ref{EVT}.
Assume, if possible, that ($i$) equivalently ($iv$) as above is false, i.e., the set of semistable points of $(G/B)^s$ is empty (and (ii)
is true). Thus,  any point $x=(\bar{g}_1,\dots,\bar{g}_s)\in (G/B)^s$
 is unstable. Choose a generic $x$. Let $\lambda=(\delta,a), P,\gamma, f, w_j$ be as above. It follows from Theorem ~\ref{jm} that $\bigcap_{j=1}^s  g_jBw_j P\subset G/P$ is the single point $f$ and, since $x$ is generic, we get
\begin{equation}\label{expdim}
[\bar{\Lambda}^{P}_{w_1}]\cdot\, \cdots \,\cdot [\bar{\Lambda}^{P}_{w_s}] = [\bar{\Lambda}^P_{e}] \in H^*(G/P,\Bbb{Z}).
\end{equation}
 We now claim that the
$s$-tuple  $(w_1, \cdots, w_s)\in (W/W_P)^s$ is $L$-movable.

Write $g_j=fp_jw_j^{-1}b_j$, for some $p_j\in P(\gamma)$ and $b_j\in B$. Hence,
$$\delta(t) \bar{g}_j=f \gamma(t) p_jw_j^{-1}B=f\gamma(t) p_j\gamma^{-1}(t)w_j^{-1}B\in G/B.$$
Define,
$$l_j = \lim_{t\to 0} \gamma(t)p_j\gamma^{-1}(t).$$
Then,  $l_j \in L(\gamma).$
Therefore,
$$\lim_{t\to 0}\delta(t)x=(f l_1 w_1^{-1}B,\dots, f l_s w_s^{-1}B).$$

By Theorem ~\ref{RR}, $\lambda\in \Lambda(\lim_{t\to 0}\delta(t)x)$.
We note that, for $j=1,\dots,s$,
$$fP(\gamma)\in (fl_j w_j^{-1})B w_j P(\gamma).$$

Applying Theorem ~\ref{jm} to the   unstable point $x_o=\lim_{t\to 0} \delta(t)x$ yields
\begin{enumerate}
\item [($\dagger$)]
 $fP$ is the only point in the  intersection
$\bigcap_{j=1}^s f l_j w_j^{-1}Bw_j P$.
\end{enumerate}
Translating by $f$, we get:
\begin{enumerate}
\item[($\ddagger$)]$eP$ is the only point in the intersection
$$\Omega:=\bigcap
l_j w_j^{-1}Bw_j P.$$
\end{enumerate}
Since the sequence  $(w_1, \cdots, w_s)$ satisfies equation (\ref{expdim}); in particular, it satisfies equation (\ref{dim0}).   Therefore, the  expected dimension of
$\Omega$ is $0$ and so is its actual dimension by ($\ddagger$).
 If this intersection $\Omega$ is not transverse at $e$,
 then by intersection theory ([Fu1, Remark 8.2]), the local multiplicity
at $e$ is $>1$, each $w_j^{-1}Bw_jP$ being  smooth.

Further, $G/P$ being a homogenous space, any other component of the intersection$$\bigcap
l_j \overline{w_j^{-1}Bw_j P}.$$
contributes nonnegatively to the intersection product
$[\bar{\Lambda}^P_{w_1}]\cdot\, \cdots \,\cdot [\bar{\Lambda}^P_{w_s}]$
(cf., [Fu1, $\S$12.2]). Thus, from equation (\ref{expdim}), we get that
the intersection   $\bigcap l_jw_j^{-1} Bw_jP$ is  transverse
 at $e\in G/P$, proving that  $(w_1,\dots,w_s)$ is $L$-movable. Clearly,
 by condition (ii), we see that the $s$-tuple $(w_1,\dots,w_s)$ satisfies
$$\sum_{j=1}^s \nu_j(w_j x_{i_P})\leq 0.$$
This  contradicts the (b)-part  of Theorem ~\ref{jm}. Thus, the set of semistable points of $(G/B)^s$ is nonempty, proving condition (i) of Theorem
\ref{EVT}.
\qed

\subsection{Proof of Theorem \ref{jm}}\label{LMIL}
Even though proof of Theorem  \ref{jm} can be extracted from [LM], we did not find it explicitly stated there. So, for completeness, we give its proof couched entirely in algebro-geometric language.

Let $x=(g_1B, \cdots, g_sB)\in (G/B)^s$ be any unstable point and let $\lambda=(\delta,a)$ be a Kempf OPS associated to $x$.
Write
$$Y=\frac{X_{\delta}}{a}=\sum_{\alpha_i\in S_\delta}c_{i} x_{i},$$
where $S_\delta$ is a subset of the set of simple roots such that  each $c_{i}$ is nonzero (and positive). Let $\gamma\in O(G)$ and $f\in G$ be such that $f\gamma f^{-1}=\delta$ and $\dot{\gamma}=X_{\delta}\in \frh_{+}$.
Define $w_j\in W/W_{P(\gamma)}$ by
$fP(\gamma)\in g_j Bw_jP(\gamma)$ for $j=1,\dots,s$.

With this notation, Theorem ~\ref{jm} can equivalently be formulated as the
 following:
\begin{enumerate}
\item[(a)$'$] For $\alpha_{i}\in S_\delta$, $fP(\lambda_{i})$ is the only point
in
$\bigcap_j g_jBw_j P(\lambda_{i})$, where $\lambda_i$ is the OPS with
$\dot{\lambda}_i=x_i$.  Moreover,
\begin{equation}\label{instabilitie}
\sum_{j=1}^s \nu_j(w_j x_{i})>0.
\end{equation}
\end{enumerate}
Note that
\begin{equation}\label{eqn7.-1}P(\lambda_{i}) \supset P(\gamma).
\end{equation}
We also note the following additional property:
\begin{enumerate}
\item[(a$_1$)] $fP(\gamma)$ is the only point  in $\bigcap g_jBw_j P(\gamma)$.
\end{enumerate}
For (a$_1$), let $hP(\gamma)$ be some other point in the intersection. Then, from Lemma \ref{l1},  we have
$\mu^{\elal}(x,h\gamma h^{-1})=\mu^{\elal}(x,f\gamma f^{-1})$. Further, clearly
$f\gamma f^{-1}$ and $h\gamma h^{-1}$ have the same norm.
Therefore, by Theorem \ref{kempf}, $fP(\gamma)f^{-1}=hP(\gamma)h^{-1}$, i.e.,
$hP(\gamma)=fP(\gamma)$ as elements of $G/P(\gamma)$.

To prove (a)$'$, it is convenient to prove the following auxiliary lemma.
\begin{lemma}\label{claim}
The inequality (\ref{instabilitie}) holds. Moreover,  let $\theta$
be  an OPS such that $X_{\theta}= x_{i}$ for $\alpha_i\in
S_\delta$ and
\begin{equation}\label{assumption}
\frac{-\mu^{\elal}(x,\theta)}{q(\dot{\theta})}\geq \frac{-\mu^{\elal}(x,f\lambda_{i}f^{-1})}{q(x_i)}.
\end{equation}
Then, equality holds in the above equation (\ref{assumption}) and $P(\theta)=P(f\lambda_{i}f^{-1})$.
\end{lemma}
This lemma will yield (a)$'$, because if $hP(\lambda_{i})$ were another
point of the intersection $\bigcap_j g_jBw_j P(\lambda_{i})$, then (by using Lemma \ref{l1} again)
$$\mu^{\elal}(x,h\lambda_{i}h^{-1})=\mu^{\elal}(x,f\lambda_{i}f^{-1})$$
and,  $q(h\lambda_ih^{-1})=q(f\lambda_if^{-1})$.
Thus, from the above lemma,  we get that $hP(\lambda_{i})=fP(\lambda_{i}).$
This completes the proof of (a)$'$.
\vskip1ex
{\bf Proof of Lemma \ref{claim}:}
Let $\tilde{H}$ be a maximal torus of $G$ contained in $P(\theta)\cap P(\delta)$, where $\lambda=(\delta,a)$ is a Kempf's OPS associated to $x$.
Using equation (\ref{translate}) and the conjugacy of maximal tori in the algebraic group $P(\theta)$, replace $\theta$ with $p\theta p^{-1}$ for some $p\in P(\theta)$ chosen so that   $p\theta p^{-1} \in O(\tilde{H})$ (and the
inequality
(\ref{assumption}) is satisfied for $\theta$ replaced by $p\theta p^{-1}$).

 Let $(\tilde{\delta},a)$ be a Kempf's OPS (corresponding to the point $x\in \xx$) such that $\tilde{\delta}\in O(\tilde{H})$. Find $b\in G$ so that  $b^{-1}\tilde{H}b=H$ and
$b^{-1}\tilde{\delta} b=X_{\tilde{\delta}}=X_{\delta}=\dot{\gamma}.$

From the uniqueness of the parabolics associated to Kempf's OPS (cf. Theorem
\ref{kempf}), we get:
\begin{equation}\label{good2}
bP(\gamma)=fP(\gamma)\in G/P(\gamma).
\end{equation}
Now, $b^{-1}\theta b\in O(H)$ and $X_{\theta}=x_i$. Therefore,
\begin{equation}\label{eqn7.0}
\Ad (b^{-1})\dot{\theta} =wx_{i}, \,\,\,\text{for some}\,  w\in W.
\end{equation}

 Let $\frh^\Bbb Q$ be the $\Bbb Q$-vector subspace of $\frh$ spanned by $\dot{\nu}$,
where $ \nu$ runs over $O(H)$, and let $\frh^\Bbb Q_+:=\frh^\Bbb Q\cap \frh_+$, where
$\frh_+:=\{h\in\frh: \alpha_i(h)\geq 0 \forall $ simple roots $\alpha_i\}$ is the set of
dominant elements of $\frh$. Define the function ${\mathfrak L}={\mathfrak
L}_{\elal,x,b}:\frh^\Bbb Q\to \Bbb{Q}$ as follows.  For any $\beta\in O(H)$ and
$r\geq 0\in \Bbb Q$, $${\mathfrak L}(r\dot{\beta})= -r\mu^\elal(x,b\beta b^{-1}).$$ Let
$V$ be a finite dimensional representation of $G$ together with a $G$-equivariant
embedding $i:\exx\to \Bbb P(V)$ such that $i^*(\mathcal O(1))$ is $G$-equivariantly
isomorphic with $\elal^M$ for some $M>0$. We can take, e.g., $V=H^0((G/B)^s, \elal^M)^*$
for any $M>0$. Define a twisted action of $G$ on $V$ via $$g\odot v=(bgb^{-1})\cdot
v,\,\,\text{for}\, g\in G, v\in V.$$
 Find a basis $\{e_1,\dots,e_n\}$ of $V$ so that, under the twisted action,
$H$ acts by the character $\eta_l$ on $e_l$, i.e., $t\odot e_l=\eta_l(t)e_l$,
for $t\in H$. Write $i(x)=[\sum x_l e_l]$. Then, by equation (\ref{muvalue}),
for any $h\in \frh^\Bbb Q$,
$${\mathfrak L}(h)=-\frac{1}{M}\max_{l:x_l\neq 0}(-\dot{\eta}_l(h))=
\frac{1}{M}\min_{l:x_l\neq 0}(\dot{\eta}_l(h)).$$
The function ${\mathfrak L}$ satisfies the following properties:
\begin{enumerate}
\item[(P$_1$)] ${\mathfrak L}$ is convex: ${\mathfrak L}(ah_1+bh_2)\geq a{\mathfrak L}(h_1)+b{\mathfrak L}(h_2)$ for $h_1,h_2\in
\frh^\Bbb Q$ and positive rational numbers $a,b$.
\item[(P$_2$)] ${\mathfrak L}(h)=\sum_j\nu_j(w_j h) $, for $h\in \frh_{+}^\Bbb Q$; in particular, ${\mathfrak L}$ is linear restricted to $\frh_{+}^\Bbb Q$.
\item[(P$_3$)] From Kempf's theory, the function  $J(h):=\frac{{\mathfrak L}(h)}{q(h)}$ on $\frh^\Bbb Q$ achieves its maximum uniquely at  the positive ray
through $Y\in \frh^\Bbb Q_+$, where $Y$ is defined in the beginning of this subsection
7.4.
\end{enumerate}
Fix $h\in \frh^\Bbb Q$ and consider the function $v \mapsto J(Y+vh),$ for
rational $v\geq 0$. Then, by the convexity of ${\mathfrak L}$ as in (P$_1$), $J(Y+vh)\geq \frac{{\mathfrak L}(Y)+v{\mathfrak L}(h)}{q(Y+vh)}$. View the right hand side as a function of $v$. It clearly takes the maximum value at $v=0$. Thus, taking its derivative
at $v=0$, we get:
$${\mathfrak L}(h)q(Y)^2-{\mathfrak L}(Y)\langle Y,h\rangle\leq 0,$$
or, that
\begin{equation}\label{a1}
J(h)\leq J(Y)\frac{\langle Y,h\rangle}{q(Y)q(h)}.
\end{equation}

From now on till the end of this proof, we take $i$ such that $\alpha_i \in S_\delta$. We  note that $Y+vx_i \in \frh_{+}$ for small
(positive or negative) values of $v$. Moreover, by (P$_2$),
${\mathfrak L} (Y+vx_i)={\mathfrak L} (Y)+v{\mathfrak L} (x_i)$ for small values of $v$. Thus, the function $v\mapsto \frac{{\mathfrak L}(Y)+v{\mathfrak L}(x_i)}{q(Y+vx_i)}$ has a local maximum at $0$; in particular, its derivative at $v=0$ is zero. This gives:
\begin{equation}\label{a2}
J(x_{i})= J(Y)\frac{\langle Y,x_i\rangle}{q(Y)q(x_i)}.
\end{equation}
Moreover,  since $\langle x_{i},x_{j}\rangle\geq 0$ for simple roots $\alpha_i$ and $\alpha_j$ and $J(Y)>0$, we get that $J(x_i)>0$. That is,
 $$\mu^{\elal}(x,b\lambda_i b^{-1})<0,$$
where $\lambda_i\in O(H)$ is defined by $\dot{\lambda_i}=x_i$.
But it follows from equation (\ref{good2}) that
\begin{equation}\label{eqn36}\mu^{\elal}(x,b\lambda_ib^{-1})=\mu^{\elal}(x,f\lambda_if^{-1}).
\end{equation}
We conclude using Lemma ~\ref{l1} that the inequality (\ref{instabilitie}) holds.
Our assumption (\ref{assumption}) now reads as
\begin{equation}\label{eqn7.1} J(wx_i)\geq J(x_i).
\end{equation} But
 according to  the inequality (\ref{a1}),
\begin{equation}\label{eqn7.2}J(wx_{i})\leq J(Y)\frac{\langle Y,wx_{i}\rangle}{q(Y)q(wx_{i})},
\end{equation}
and by equation (\ref{a2}),
   \begin{equation}\label{eqn7.3}J(x_{i})= J(Y)\frac{\langle Y,x_{i}\rangle}{q(Y)q(x_{i})}.
\end{equation}

It is also easy to see that
$\langle Y,x_{i}\rangle >\langle Y,wx_{i}\rangle$, if $wx_{i}\neq x_{i}$. Combining equations (\ref{eqn7.1})--(\ref{eqn7.3}),  we therefore conclude that $wx_i=x_i$. Thus, by equation   (\ref{eqn7.0}), $\theta=b\lambda_ib^{-1}$ and hence
$$P(\theta)=P(b\lambda_ib^{-1})=P(f\lambda_if^{-1}),$$
by equations (\ref{eqn7.-1}) and (\ref{good2}). Finally, by
equations (\ref{eqn7.0}) and (\ref{eqn36}), the inequality (\ref{assumption})
is in fact an equality. This proves  Lemma \ref{claim}.
\qed

\subsection{$L$-movability and the eigenvalue problem}

Let $G$ be a connected semisimple group. Choose  a maximal compact subgroup
  $K$ of $G$ with Lie algebra $\frk$. Then,
there is a natural homeomorphism
$C:\frk/K\to \frh_{+}$, where $K$ acts on $\frk$ by the adjoint representation.

Now, one of the main aims of the {\it eigenvalue problem} is to describe the set $\Gamma(s,K):=$
$$\{(h_1,\dots,h_s)\in\frh_{+}^s\mid \exists (k_1,\dots,k_s)\in \mathfrak k^s \text{:} \sum_{j=1}^s k_j=0\,\,\text{and }\, C(k_j)=h_j \forall j=1,\dots,s\}.$$

Given a standard maximal parabolic subgroup $P$, let $\omega_P$ denote the corresponding fundamental weight, i.e., $\omega_P(\alpha_i^\vee)=1$, if $\alpha_i \in \Delta\setminus \Delta(P)$ and $0$ otherwise, where $\alpha_i^\vee$ is the
fundamental coroot corresponding to the simple root $\alpha_i$. This is invariant under the Weyl group $W_P$ of $P$.

The following theorem is one of our main results which gives a solution of the
eigenvalue problem.
\begin{theorem}\label{eigen}  Let $(h_1,\dots,h_s)\in\frh_{+}^s$.  Then, the following
are equivalent:

(a) $(h_1,\dots,h_s)\in\Gamma(s,K)$.

(b)
 For every standard maximal parabolic subgroup $P$ in $G$ and every choice of
  $s$-tuples  $(w_1, \cdots, w_s)\in (W^P)^s$ such that
$$[\bar{\Lambda}^P_{w_1}]\odot_0\, \cdots \,\odot_0 [\bar{\Lambda}^P_{w_s}] = [\bar{\Lambda}^P_{e}] \in \bigl(H^*(G/P,\Bbb{Z}), \odot_0\bigr),$$
the following inequality holds:
$$\omega_P(\sum_{j=1}^s\,w_j^{-1}h_j)\leq 0.$$
\end{theorem}

\begin{proof} Observe first that, under the identification of $\frh$ with $\frh^*$ induced from the Killing form,   $\frh_+$ is isomorphic with  the set
$D$ of dominant weights of $\frh^*$. In fact, under this identification, $x_i$
corresponds with  $2\omega_i/\langle  \alpha_i, \alpha_i\rangle$, where $\omega_i$ denotes the $i$-th fundamental weight. Let $D_\Bbb Z$ be the set of
dominant integral weights. Define
\begin{align*} \bar{\Gamma}(s)&:=\{(\nu_1, \cdots, \nu_s)\in D^s: N\nu_j\in D_\Bbb Z
\,\,\text{for all}\, j\,\,\text{and}\\
& H^0\bigl((G/B)^s,
\cl(N\nu_1)\boxtimes \cdots \boxtimes\cl(N\nu_s)\bigr)^G\neq 0\,\text {for
some}\, N>0\}.
\end{align*}

Then, under the identification of $\frh_+$ with $D$ (and hence of $\frh_+^s$ with $D^s$),
$\Gamma (s,K)$ corresponds to the closure of $\bar{\Gamma}(s)$. In fact,
$\bar{\Gamma}(s)$ consists of  the rational points of the image of $\Gamma(s,K)$
(cf., e.g., [Sj, Theorem 7.6]). Since $x_i$ corresponds with $2\omega_i/\langle  \alpha_i, \alpha_i\rangle$, the theorem follows from Theorem \ref{EVT}.

\end{proof}

\section{Nonvanishing of Products in the  Cohomology of Flag Varieties (Horn Inequalities)}

We give two inductive criteria (actually, only necessary conditions) to determine when the product of a number of Schubert cohomology classes of $G/P$ is nonzero. The first criterion (Theorem \ref{T2}) is in terms of the characters, whereas the second one (Theorem \ref{weakHorn}) is in terms of dimension counts.

\begin{theorem}\label{T2}
Assume that  $(w_1, \cdots, w_s)\in (W^P)^s$  satisfies equation (\ref{dim0}) and that
$[\bar{\Lambda}^P_{w_1}]\cdot\, \cdots \,\cdot [\bar{\Lambda}^P_{w_s}] = d[\bar{\Lambda}^P_{1}]$ in  $H^*(G/P)$, for some nonzero $d$.
Then,

(a)  For each $\alpha_i\in \Delta\setminus\Delta(P)$, the following inequality holds:
\begin{equation}\label{center}
\bigl((\sum_{j=1}^s\,\chi_{w_j})-\chi_1\bigr)(x_{i})\leq 0,
\end{equation}
where $\chi_w$ is defined in Definition \ref{basicbundle}.

(b)  For any standard  parabolic $Q_L$ of $L$ (i.e., $Q_L\supset B_L$), and
 $u_1, \cdots, u_s\in W_L/W_{Q_L}$ such that
$$
[\bar{\Lambda}^{Q_L}_{u_1}]\cdot \,\cdots \,\cdot [\bar{\Lambda}^{Q_L}_{u_s}]
\neq 0\in H^*(L/Q_L),$$
the inequality
\begin{equation}\label{center'} \sum_{j=1}^s \,\chi_{w_j}(u_jx_p)\leq
\chi_1(x_p)\end{equation}
holds for any $p$ such that
  $\alpha_p \in \Delta(P)\setminus\Delta(Q_L)$.
\end{theorem}

\begin{proof}
Let  $\nu_1, \cdots, \nu_s \in X(H)$ and let
$$0\neq\sigma\in H^0\left( (P/B_L)^{s}, \cl (\nu_1) \boxtimes \cdots \boxtimes \cl (\nu_s)\right)^P.$$
Assume  that $\sigma$ does not vanish at $(\bar{p}_1, \cdots, \bar{p}_s)\in (P/B_L)^s$.
Then, for every admissible OPS $\lambda \in O(P)$, we have the following inequality obtained from Proposition
\ref{propn14}(b) and Lemma \ref{l1}:
\begin{equation}\label{ES}
\nu_1([\bar{p}_1,\lambda]{X}_{\lambda})+\cdots
+\nu_s([\bar{p}_s,\lambda]{X}_{\lambda})\leq 0,
\end{equation}
where ${X}_\lambda$ is defined by equation (\ref{Xtilde}).

 Now, Theorem ~\ref{T2} follows immediately from the following proposition
together with Proposition \ref{First} and Lemma \ref{sectiontheta}. (Observe that $\chi_1$ is fixed by $W_L$.)

\end{proof}
\begin{proposition}\label{nonzerosections}
Let $\nu_1, \cdots, \nu_s \in X(H)$ and let
\begin{equation}\label{vs}
H^0\left( (P/B_L)^{s}, \cl (\nu_1) \boxtimes \cdots \boxtimes \cl (\nu_s)\right)^P\neq 0.
\end{equation}
 Then,

(a) For each $\alpha_i\in \Delta\setminus\Delta(P)$, we have
$\sum^s_{j=1} \nu_j(x_i) \leq 0$.

(b) For any  standard  parabolic subgroup $Q_L$ of $L$ and $u_1, \cdots, u_s \in W_L/W_{Q_L}$ such that $[\bar{\Lambda}^{Q_L}_{u_1}]\cdot \,
\cdots  \,\cdot[\bar{\Lambda}^{Q_L}_{u_s}]\neq 0\in H^*(L/Q_L)$, we have:
$$\sum_{j=1}^s \nu_j(u_j x_{p})\leq 0, \,\,\text{for all}\,\, \alpha_p \in \Delta(P)\setminus\Delta(Q_L).$$
\end{proposition}
\begin{proof}
For (a), apply equation (\ref{ES}) to the admissible OPS
$\lambda=t^{x_i}$, which is
central in $L$. In this case, by equation (\ref{eqn11}), $[\bar{p},\lambda]=1$
for any $\bar{p}\in P/B_L$  and $X_\lambda=\dot\lambda=x_i$.

To prove (b), pick a nonzero $\sigma$ in the vector space (\ref{vs}). Let
$Z\subset L^s$ be the set of  points $(l_1, \cdots, l_s)$ such that
$$l_1B_L u_1 Q_L \cap \cdots \cap l_sB_L u_s Q_L \neq \emptyset.$$
By the assumption in (b), $Z$ is nonempty (and open).  Let $Z_{\sigma}$ be the subset of $P^s$ consisting of  $(p_1, \cdots, p_s)$ such that $\sigma$
does not vanish at  $(\bar{p}_1, \cdots, \bar{p}_s)\in (P/B_L)^s$, where $\bar{p}_j:=p_jB_L$. Since $\sigma \neq 0$, $Z_\sigma$ is a
nonempty open subset of $P^s$. Consider the projection $\pi:P^s\to L^s$
under the decomposition $P=U\cdot L$  and pick
$$(p_1, \cdots, p_s)\in Z_\sigma\cap\pi^{-1}(Z).$$
Since $Z$ and $Z_\sigma$ are nonempty Zariski open subsets, the intersection
$Z_\sigma\cap\pi^{-1}(Z)$ is  nonempty and open.

Let $p_j=u_j l_j,$ for  $j=1, \cdots, s$, where $u_j\in U$ and $l_j\in L$.  Pick
$$l\in l_1B_L u_1 Q_L \cap \cdots \cap l_sB_L u_s Q_L. $$
Consider the admissible OPS $\lambda (t)=lt^{x_{p}}l^{-1},$  for $\alpha_p \in \Delta(P)\setminus\Delta(Q_L).$ (To prove that this is
admissible, use Lemma \ref{lemmar}.) Since $\lambda$  is conjugate to the OPS $\lambda_o:=t^{x_p}$ lying in $H$ and, moreover, $x_p\in \fh$ is $L$-dominant,
we get $X_\lambda=x_p$. Clearly,
$[\bar{p}_j, \lambda]=u_j$,  for any  $j=1,\cdots, s$. We can therefore use equation
(\ref{ES}) to conclude that
$$\sum_{j=1}^s\nu_j( u_jx_{p})\leq 0.$$
This proves (b).

\end{proof}
\begin{remark} (a)  Let $\nu_1, \cdots, \nu_s \in X(H)$ and let $\elal$ be the line bundle  $\cl (\nu_1) \boxtimes \cdots \boxtimes \cl (\nu_s)$ on $(P/B_L)^{s}$. Assume that for every $\alpha_i\in \Delta\setminus\Delta(P)$, we have
$\sum^s_{j=1} \nu_j(x_i) =0$. Then, the restriction map 
$$H^0\left( (P/B_L)^{s}, \elal\right)^P \to H^0\left(
(L/B_L)^{s},  \elal\right)^L$$
is an isomorphism. If the above equality is violated for some 
 $\alpha_i\in \Delta\setminus\Delta(P)$, then  $H^0\left(
(L/B_L)^{s},  \elal\right)^L=0$.

 This is proved by
the same technique that was used in the proof of Theorem ~\ref{T1} (cf., Section
~\ref{ref11}).
 
(b)  In Theorem \ref{T2}, the validity of property (b) for every standard parabolic 
subgroup 
$Q_L$ of $L$ is equivalent to the corresponding property only for the standard {\it maximal} parabolic subgroups $Q_L$ of $L$. This can be proved by the same technique as developed in Section 7.
\end{remark}

In the $L$-movable case, we have the following refinement of Theorem \ref{T2}. To state this refinement, we need the following notation.

For any $w\in W^P$ and a central character $c$ of $L$ (i.e., an algebraic group homomorphism $ Z(L)\to \Bbb G_m$, $Z(L)\subset H$ being the center of the Levi subgroup $L$ of $P$), define
\begin{equation}\label{eqnm1} \chi_w^c=\sum_{\beta\in R(w,c)}\,\beta,
\end{equation}
where
\begin{equation}\label{eqnm2}  R(w,c):=\{\beta \in \bigl(R^+
\setminus R^+_{\fl}\bigr)
\cap w^{-1}R^+: e^\beta_{\mid Z(L)} =c\}.
\end{equation}
Observe that
\begin{equation}\label{eqnm3} \chi_w=\sum_c\chi_w^c,
\end{equation}
where the sum runs over all the central characters of $L$ such that
$\chi_1^c\neq 0$.

\begin{theorem}\label{T2'}
Assume that  the $s$-tuple $(w_1, \cdots, w_s)\in (W^P)^s$  is
$L$-movable. Then,
\begin{enumerate}
\item For any central character $c$ of $L$ such that $\chi_1^c\neq 0$, we have
\begin{equation}\label{addition} \sum_{j=1}^s \,\mid R(w_j,c)\mid=
\mid R(1,c)\mid,\end{equation}
where $\mid \cdot\mid$ denotes the cardinality of the enclosed set.
\item For any standard  parabolic $Q_L$ of $L$  and
 $u_1, \cdots, u_s\in W_L/W_{Q_L}$ such that
$$
[\bar{\Lambda}^{Q_L}_{u_1}]\cdot \,\cdots \,\cdot [\bar{\Lambda}^{Q_L}_{u_s}]
\neq 0\in H^*(L/Q_L),$$
and any central character $c$ of $L$ such that $\chi_1^c\neq 0$, the following
 inequality is satisfied for any
  $\alpha_p \in \Delta(P)\setminus\Delta(Q_L)$:
\begin{equation}\label{eqnm4} \sum_{j=1}^s \,\chi_{w_j}^c(u_jx_p)\leq
\chi_1^c(x_p).\end{equation}
\end{enumerate}
\end{theorem}
\begin{remark}
(a) Observe that in the $L$-movable case, by virtue of Theorem \ref{T1}, the inequality (\ref{center})
is, in fact, an equality.

(b) The inequalities (\ref{eqnm4}) summed over all the central characters $c$ of $L$ such that $\chi_1^c\neq 0$ is nothing but the inequality (\ref{center'})
(use the identity (\ref{eqnm3})). Thus, Theorem \ref{T2'} is a refinement of Theorem \ref{T2}  in the $L$-movable case.
\end{remark}

\begin{proof} (of Theorem \ref{T2'}) For any  central character $c$ of $L$
and $w\in W^P$, let $T^c_w$ be the $B_L$-submodule of $T_w$ defined by
$$ T_w^c:=\{v\in T_w: t\cdot v=c(t)v, \forall\text{ }t\in Z(L)\subset B_L\}.$$
This gives rise to the $L$-equivariant vector bundle on $L/B_L$:
$$\mathcal{T}_w^c:=L\times_{B_L}\,T_w^c.$$
Similarly, we can define the $B_L$-submodule  $T^c$ of $T$
and the associated  $L$-equivariant vector bundle $\mathcal{T}^c$ on $L/B_L$
and $\mathcal{T}^c_s$ on $(L/B_L)^s$.
Analogous to Lemma \ref{mt}, we have:
\begin{equation}\label{eqnm5} \det (\mt^c/\mathcal{T}_w^c)=\cl (\chi_w^c),
\end{equation}
as $L$-equivariant vector bundles on $L/B_L$. Let $\Theta_o$
denote the restriction of the bundle map $\Theta$ (defined by
(\ref{map})) to the subvariety $(L/B_L)^s\subset (P/B_L)^s$. Then,
by Corollary \ref{movL}, $\Theta_o$ is an isomorphism over a dense
open subset of $(L/B_L)^s$. From this we see that, for any central
character $c$ of $L$ such that $\chi_1^c\neq 0$,
\begin{equation}\label{eqnm6}
{\Theta_o}_{\mid \mathcal{T}^c_s}:\mt_s^c \to \oplus_{j=1}^s
\pi_j^*({\mt}^c/\mathcal{T}_{w_j}^c)
\end{equation}
is an isomorphism over a dense open subset of $(L/B_L)^s$,
 where   $\pi_j: (L/B_L)^s  \to L/B_L$ is the projection onto the $j$-th
factor. Therefore, the ranks of the two sides of equation (\ref{eqnm6}) coincide. This 
gives the equality (\ref{addition}).

Taking the determinant of ${\Theta_o}_{\mid \mathcal{T}^c_s}$ and using the equation
(\ref{eqnm5}), we get a nonzero bundle map
$$\theta_o^c: \det(\mt_s^c) \to
\cl (\chi^c_{w_1}) \boxtimes \cdots \boxtimes \cl (\chi^c_{w_s}),$$
i.e., a nonzero section of the line bundle
$$\cl (\chi^c_{w_1}-\chi^c_1) \boxtimes \cdots \boxtimes \cl (\chi^c_{w_s})$$
on $(L/B_L)^s$. Now, applying Proposition \ref{nonzerosections} for the case
$G=L, P=L$, and observing that $\chi_1^c$ is fixed by $W_L$, we get
the theorem.
\end{proof}

   \begin{question}  We would like to ask if the converse of Theorem \ref{T2} is true. Specifically, take  $(w_1, \cdots, w_s)\in (W^P)^s$ satisfying the equation
(\ref{dim0}) and assume that the conditions (a) and (b) of Theorem \ref{T2} are satisfied. Then, is it true that
$$[\bar{\Lambda}^P_{w_1}]\cdot\, \cdots \,\cdot [\bar{\Lambda}^P_{w_s}] = d[\bar{\Lambda}^P_{1}]$$
 in  $H^*(G/P)$, for some nonzero $d$. For $G=SL_n$ and any maximal parabolic subgroup $P$, this question has an affirmative answer [Fu2].

One could ask the corresponding (weaker) question for $L$-movable $s$-tuples, i.e., is the converse of Theorem \ref{T2'} true? It may be remarked here that this question (in the $L$-movable case) for $G/B$ has an affirmative answer by virtue of Corollary \ref{leviprod}.
\end{question}

We now come to our second criterion (Theorem \ref{weakHorn}) to determine when the products of cohomology classes in $G/P$ are nonzero. This criterion is
in terms of inequalities involving certain dimension counts.

Let $Q$ be a standard  parabolic subgroup of $G$ contained in $P$ and let
$Q_L:=Q\cap L$ be the associated  parabolic in $L$. Let $\haq$ be a standard parabolic in $G$ containing $Q$. Then, we have the canonical identification
 $L/Q_L\simeq P/Q$ and the standard projection
$\tau: P/Q\to G/\hat{Q}$.
\begin{lemma}\label{Bruhaha}
For any  $w\in W^P, v\in W_L,  p\in P$  and $b\in B$, we have (setting
 $g^{-1}= bwp$):
$$\tau(p^{-1}B v Q) \subset gB(wv)\haq .$$
\end{lemma}
\begin{proof} Since $w\in W^P$, by equation (\ref{eqn1}), we have $wB_Lw^{-1}\subset B$. Further,
$$\tau(p^{-1}B v Q)= \tau(p^{-1}B_LU v Q)=\tau(p^{-1}B_L v Q)=gbwB_Lv \haq=gbwB_Lw^{-1}(wv)\haq.$$
 We therefore have
$$\tau(p^{-1}B v Q)\subset gbBwv\haq = gB(wv)\haq.$$

\end{proof}
\begin{theorem}\label{weakHorn} Let $w_1,\dots,w_s \in W^P$ be such that
\begin{equation}\label{eqn9.1}[\bar{\Lambda}_{w_1}^P]\cdot \,\cdots \,\cdot [\bar{\Lambda}_{w_s}^P]\neq 0
\end{equation}
in $H^{*}(G/P,\Bbb{Z})$. Then given the data:
$Q\subset P$ a  parabolic, $\hat{Q}\subset G$ a  parabolic containing $Q$, elements $u_1,\dots,u_s\in W_L/W_{Q_L}$ such that
\begin{equation}\label{eqn9.2}[\bar{\Lambda}^{Q_L}_{u_1}]\cdot\, \cdots \,\cdot [\bar{\Lambda}^{Q_L}_{u_s}]\neq 0
\end{equation}
in $H^{*}(L/Q_L,\Bbb{Z})$, the following hold (with  $\hat{w}_j:=w_ju_j$):
\vskip1ex

(a)  $[\bar{\Lambda}^{\hat{Q}}_{\hat{w}_1}]\cdot \,\cdots \,\cdot
[\bar{\Lambda}^{\hat{Q}}_{\hat{w}_s}]\neq 0 $ in $H^{*}(G/\hat{Q},\Bbb{Z})$.
In particular,
\begin{equation}\label{eqn9.3}\sum_j \codim(\Lambda^{\hat{Q}}_{\hat{w}_j};G/\hat{Q})\leq \dim(G/\hat{Q}).
\end{equation}

(b) If $\hat{Q}\cap P=Q$ then,
\begin{equation}\label{eqn9.4}\dim(G/\hat{Q})-\sum_j \codim(\Lambda^{\hat{Q}}_{\hat{w}_j};G/\hat{Q})\geq \dim(L/Q_L)-\sum_j
\codim(\Lambda^{Q_L}_{{u}_j};L/Q_L).\end{equation}

The above inequality (\ref{eqn9.4}) is equivalent to the following inequality
(in the case $\hat{Q}\cap P=Q$):
\begin{equation}\label{eqn9.5}
\mid R(\fru_{\hat{Q}})\cap R(\fru_{P})\mid\geq \sum_{j=1}^s\mid R(\fru_{\hat{Q}})\cap R(\fru_{P})\cap \hat{w}_j^{-1} R^+\mid ,
\end{equation}
where $\mid \cdot\mid$ denotes the cardinality of the enclosed set.
 \end{theorem}
\begin{proof} Take generic $g_j$ for $j=1,\dots,s$ so that for each standard parabolic $\tilde{P}$ in $G$ and any $(z_1,\dots, z_s)\in W^s$, the
intersection
$${g_1Bz_1\tilde{P}}\cap \dots\cap {g_sBz_s\tilde{P}}$$ is
transverse (possibly empty).

We may further assume, by left multiplying each $g_j$ by the same element, that for the given $s$-tuple   $(w_1,\dots,w_s) \in (W^P)^s$,
$$e\in g_1Bw_1P\cap \dots \cap g_sBw_sP.$$

Choose  $p_j\in P, b_j\in B$  such that $e=g_jb_jw_jp_j$, for $j=1,\dots,s$.
Now, consider  the intersection in $L/Q_L\simeq P/Q$:
$$\Omega= \overline{p_1^{-1}Bu_1 Q}\cap \dots\cap \overline{p_s^{-1}Bu_s Q}.$$
Then, $\Omega$ is nonempty because of the  assumption (\ref{eqn9.2}). Each irreducible component of $\Omega$ is of dimension at least
\begin{equation}\label{EQ2}
\dim(L/Q_L)-\sum_{j=1}^s \codim(\Lambda^{Q_L}_{{u}_j};L/Q_L).
\end{equation}

By Lemma ~\ref{Bruhaha}, $\tau(p_j^{-1}Bu_j Q)\subset  g_jB w_ju_j\haq
\subset \overline{g_jB \hat{w}_j\haq}$, for $j=1,\dots,s$,
under the projection $\tau:P/Q \to G/\hat{Q}$ (which is an embedding if $\hat{Q}\cap P=Q$). So,
\begin{equation}\label{eqn45}\tau(\Omega)\subset \bigcap_{j=1}^s\overline{g_jB \hat{w}_j\haq}.\end{equation}
The intersection $\Omega':=\bigcap_{j=1}^s\overline{g_jB \hat{w}_j\haq}$ is therefore nonempty (and proper by the choice of $g_j$). This gives the part
(a) of the theorem.

In the case $\haq\cap P=Q$, $\tau$ is an embedding. Comparison of equations (\ref{EQ2}) and (\ref{eqn45})    gives the inequality (\ref{eqn9.4}) of (b). The inequality (\ref{eqn9.5}) follows from (\ref{eqn9.4}) and the following Lemma \ref{dimension}.
\end{proof}
\begin{remark}
This theorem is of interest  even in the case $P=B$. In this case,
$Q=B$ and any standard parabolic is allowed for $\hat{Q}$.
\end{remark}

We recall  the following simple lemma.
\begin{lemma}\label{lemma34}
For $w\in W$ (not necessarily in $W^P$), we have:
\begin{equation}\label{eqn9.6}
\codim(\Lambda_w^P;G/P)=\mid R(\fru_P)\cap (w^{-1}R^+)\mid.
\end{equation}
\end{lemma}

\begin{lemma}\label{dimension}
Let $w\in W^P$, $u\in W_P$ and $\hat{w}:=wu$. Let $Q\subset P$ be a standard
parabolic and let $\haq$ be a parabolic in $G$ such that $\haq\cap P=Q$. Then,
\begin{equation} \label{eqn9.7}
\codim(\Lambda_{\hat{w}}^{\haq}; G/\haq)-
\codim(\Lambda_u^{Q_L};L/Q_L)=
\mid R(\fru_{\hat{Q}})\cap R(\fru_{P})\cap \hat{w}^{-1} R^+\mid.
\end{equation}
\end{lemma}
\begin{proof}
Write $R^{+}$ as a disjoint union of three parts (note that since $w\in W^P$, $wR_\frl^+\subset R^+$ by (\ref{eqn1}))
$$R^+=wR_\frl^{+}\sqcup \bigl(wR(\fru_P^-)\cap R^+\bigr)\sqcup
 \bigl(wR(\fru_P)\cap R^+\bigr),
$$ where $\fru_P^-$ is the opposite nil-radical of $\frp$.
Take $\hat{w}^{-1}=u^{-1}w^{-1}$ of this decomposition and intersect with $R(\fru_{\hat{Q}})$. The first piece is
$$R(\fru_{\hat{Q}})\cap u^{-1}R_\frl^+=R(\fru_Q)\cap u^{-1}R_\frl^+,\,\,
\text{since}\, \haq\cap P=Q.$$
Thus, by Lemma \ref{lemma34},
\begin{equation}\label{eqn9.8}
\mid R(\fru_{\hat{Q}})\cap u^{-1}R_\frl^+\mid = \codim(\Lambda_u^{Q_L};L/Q_L).
\end{equation}
The second piece is
$$R(\fru_{\hat{Q}})\cap u^{-1}R(\fru_P^-)\cap u^{-1}w^{-1}R^+.$$
But $u^{-1}R(\fru_P^-)=R(\fru_P^-)$  and clearly
\begin{equation}\label{eqn9.9}R(\fru_{\hat{Q}})\cap R(\fru_P^-)=\emptyset .
\end{equation}
So the second piece gives us the empty set.

The third piece is
\begin{equation}\label{eqn9.10}R(\fru_{\hat{Q}})\cap u^{-1}R(\fru_P)\cap u^{-1}w^{-1}R^+=
R(\fru_{\hat{Q}})\cap R(\fru_P)\cap u^{-1}w^{-1}R^+.
\end{equation}
Finally, by Lemma \ref{lemma34},
\begin{equation}\label{eqn9.11}
\mid R(\fru_{\hat{Q}})\cap u^{-1}w^{-1}R^+\mid=
\codim(\Lambda_{\hat{w}}^{\haq}; G/\haq).
\end{equation}
Combining the equations (\ref{eqn9.8})--(\ref{eqn9.11}), we get equation (\ref{eqn9.7}), proving the lemma.
\end{proof}
\begin{remark}
(a) It is easy to see that in the case when $\haq$ is a minuscule maximal parabolic subgroup and $Q\neq P$ (and $\haq \cap P=Q$), the inequality
(\ref{eqn9.5}) is the same as the inequality (\ref{center'}).

Also, in the case $\haq=Q$ and $P$ is a minuscule maximal parabolic subgroup,
the inequality (\ref{eqn9.5}) is the same as the inequality (\ref{center}).

It may be remarked that, in general, the inequality (\ref{center'}) is a 
certain
`weighted' version of  the inequality (\ref{eqn9.5}).

(b) In the general case, we do not know if the system of inequalities
(\ref{center}) and (\ref{center'}) is equivalent to the system of inequalities
(\ref{eqn9.5}). In fact, we would expect that the former set is more refined than the latter.

(c) In his PhD thesis ``Vanishing and non-vanishing criteria for branching Schubert calculus,''  Kevin Purbhoo has given some criteria for determining which of the Schubert
intersection numbers vanish in terms of a combinatorial recipe which he calls `root game.'
\end{remark}

 \section{A further study of $ \bigl( H^*(G/P, \bc ), \odot_0\bigr)$}

For any Lie algebra $\fs$ and a subalgebra $\ft$, let $H^*(\fs ,\ft )$
be the Lie algebra cohomology of the pair $(\fs ,\ft )$ with
trivial coefficients.  Recall (cf. [Ku2, \S 3.1]) that this is the cohomology of the cochain
complex
  \begin{align*}
C\u. (\fs ,\ft ) &= \{ C^p(\fs ,\ft )\}_{p\geq 0}, \quad\text{ where}\\
C^p(\fs ,\ft ) &:= \Hom_{\ft}\bigl( \wedge^p(\fs /\ft ), \bc \bigr).
  \end{align*}

We now return to the notation of \S 2.  For any (positive) root $\beta\in
R^+$, let $y_{\beta}\in \fg_{\beta}$ be the corresponding Chevalley
root vector  (which is unique up to a sign) and let $y_{-\beta}\in\fg _{-\beta}$ be
the vector such that $\ip<y_{\beta}, y_{-\beta}> =1$ under the Killing
form. For any $w\in W^P$, let $\Phi_w := w^{-1}R^-\cap R^+ \subset R(\fu_P)$.  Then,
\begin{equation}\sum_{\beta\in\Phi_w} \beta = \rho -w^{-1}\rho.
\end{equation}
In particular,  $\Phi_v =
\Phi_w$ iff $v=w$.
Let $\Phi_w = \{ \beta_1,\cdots ,\beta_p\} \subset
R(\fu_P)$.  Set $y_w := y_{\beta_1}\wedge \cdots\wedge y_{\beta_p}\in
\wedge^p(\fu_P)$, determined up to a sign.  Then, up to scalar multiples,
$y_w$ is the unique weight vector of $\wedge (\fu_P)$ with weight $\rho
-w^{-1}\rho$.  Similarly, we can define $y^-_w\in\wedge^p(\fu^-_P)$ of
weight $w^{-1}\rho -\rho$.

\vskip1ex
We recall the following fundamental result due to Kostant  [Ko1].

\begin{theorem}\label{Th9.A}  For any standard parabolic subgroup $P$ of
$G$,
  \[
H^p (\fu_P) = \bigoplus_{\substack{ w\in W^P:\\ \ell (w)=p}} M_w,
  \]
 as $\fl$-modules, where $M_w$ is the unique irreducible $\fl$-submodule of $H^p(\fu_P)$ with
highest weight $w^{-1}\rho -\rho$ (which is $\fl$-dominant for any
$w\in W^P$). This has  a highest weight vector $\phi_w\in\wedge^p(\fu_P)^*$ defined
by $\phi_w(y_w)= 1$ and $\phi_w(y)=0$ for any weight
vector of $\wedge^p(\fu_P)$  of weight $\neq \rho-w^{-1}\rho$.

Similarly, for the opposite nil-radical $\fu^-_P$,
  \[
H^p(\fu^-_P) = \bigoplus_{\substack{ w\in W^P:\\ \ell (w)=p}} N_w,
  \]
 as $\fl$-modules, where $N_w$ is the unique  irreducible $\fl$-submodule of
$H^p(\fu_P^-)$ isomorphic
with the dual $M^*_w$ and it has a lowest weight vector
$\phi^-_w\in\wedge^p(\fu^-_P)^*$ defined by $\phi_w^-(y_w^-)= 1$ and $\phi_w^-(y)=0$ for any weight
vector of $\wedge^p(\fu_P^-)$  of weight $\neq w^{-1}\rho -\rho$.

Thus,
  \begin{align*}
[H^p(\fu_P)\otimes H^q(\fu^-_P)]^{\fl} &=0, \quad\text{unless $p=q$, and}\\
[H^p(\fu_P)\otimes H^p(\fu^-_P)]^{\fl} &\simeq \bigoplus_{\substack{ w\in
W^P,\\ \ell (w)=p}} \Bbb C\xi^w,
  \end{align*}
where $\xi^w\in [M_w\otimes N_w]^{\fl}$ is the unique element whose
$H$-equivariant projection to $(M_w)_{w^{-1}\rho -\rho} \otimes N_w$ is the
element $\phi_w\otimes\phi_w^-$, $(M_w)_{w^{-1}\rho -\rho}$ being the weight
space of $M_w$ corresponding to the weight $w^{-1}\rho -\rho$.  (Observe that
the ambiguity in the sign of $y_w$ disappears in the definition of
$\xi^w$ giving rise to a completely unique element.)
\end{theorem}

  \begin{theorem}\label{liecohomology}  For any standard parabolic subgroup $P$ of $G$, there
is a graded algebra isomorphism
 $$\phi : \bigl( H^*(G/P, \bc ), \odot_0\bigr) \simeq \bigl[
H^*(\fu_P)\otimes H^*(\fu_P^-)\bigr]^{\fl}$$
such that
 \begin{equation}\label{eqn10.0} \phi\bigl([\bar{\Lambda}^P_w]\bigr) = \Bigl( \frac{i}{2\pi}\Bigr)^{\dim
G/P-\ell (w)} \xi^{w_oww_o^P} \prod_{\al\in\Phi_{w_oww_0^P}} \ip<\rho
,\al>,
\end{equation}
where $H^*(G/P,\bc )$ is equipped with the product $\odot_0$
as in Definition ~\ref{levicohomology}, and we take the tensor product algebra structure on
the right side.
  \end{theorem}

  \begin{proof}
Let $x_P := \sum_{\al_i\in\Del\backslash\Del (P)}x_i$.  Consider the
graded, multiplicative decreasing filtration $\ca =\{\ca_m\}_{m\geq 0}$ of
$H^*(G/P, \bc )$ under the intersection product defined as follows:
  \[
\ca_m := \bigoplus_{w\in W_m^P} \bc[\bar{\Lambda}^P_w],
  \]
where $W_m^P := \{ w\in W^P: \chi_w(x_P)\geq m\}$.  Observe that, by equation
(\ref{eqn5}), $\chi_w(x_P)= \rho +w^{-1}\rho (x_P)\in\bz_+$.  Thus, by Proposition
18(a), $\ca_m\ca_n\subset \ca_{m+n}$, i.e., $\ca$ respects the algebra
structure.

Let $\gr(\ca ) := \bigoplus_{m\geq 0}\frac{\ca_m}{\ca_{m+1}}$ be the
associated `gr' algebra.  Then $\gr(\ca )$ acquires two gradings: the
first one the `homogeneous grading' $m$ assigned to the elements of
$\frac{\ca_m}{\ca_{m+1}}$ and the second one the `cohomological grading'
coming from the cohomological degree in $H^*(G/P, \bc )$.  For example,
$[\bar{\Lambda}_w^P]$ has bidegree $(\chi_w(x_P), \dim G/P-\ell (w))$.

By the definition of the product $\odot_0$ and Proposition 18(a), the
linear map
  \[
\varphi : \bigl( H^*(G/P,\bc ), \odot_0\bigr) \to \gr(\ca ),
[\bar{\Lambda}^P_w] \mapsto [\bar{\Lambda}^P_w] \mod \ca_{\chi_w(x_P)+1},
  \]
is, in fact, a graded algebra isomorphism with respect to the
cohomological grading on $\gr(\ca )$.

We next introduce another filtration $\{\bar\cf_m\}_{m\geq 0}$ of
$H^*(G/P,\bc )$ in terms of the Lie algebra cohomology.  Recall
that choosing a maximal compact subgroup $K$ of $G$, we can
identify $H^*(G/P,\bc )$ with the cohomology of the $K$-invariant
forms on $G/P$ under the de Rham differential, i.e., with the Lie
algebra cohomology $H^*(\fg ,\fl )$.  The underlying cochain
complex $C\u. =C\u. (\fg ,\fl )$ for $H^*(\fg ,\fl )$ can be
rewritten as
  \[
C\u. := [\wed\u. (\fg /\fl )^*]^{\fl} = \Hom_{\fl} \bigl( \wed\u.
(\fu_P)\otimes\wed\u. (\fu_P^-), \bc \bigr).
  \]
Define a decreasing filtration $\cf =\{\cf_m\}_{m\geq 0}$ of the
cochain
complex $C\u.$ by subcomplexes:
  \[
\cf_m := \Hom_{\fl}\Biggl( \frac{\wed\u. (\fu_P)\otimes\wed\u.
(\fu^-_P)}{\bigoplus_{s+t\leq m-1}
\wed\u._{(s)}(\fu_P)\otimes\wed\u._{(t)}(\fu^-_P)}, \bc\Biggr) ,
  \]
where $\wed\u._{(s)}(\fu_P)$ (resp. $\wed\u._{(s)}(\fu^-_P)$) denotes the
subspace of $\wed\u. (\fu_P)$ (resp. $\wed\u. (\fu^-_P)$) spanned by the
$\fh$-weight vectors of weight $\beta$ with {\em $P$-relative height}
$$\text{ht}_P (\beta ) := \mid\beta (x_P)\mid=s.$$

Now, define the filtration $\bar\cf =\{ \bar\cf_m\}_{m\geq 0}$ of
$H^*(\fg ,\fl )\simeq H^*(G/P)$ by
  \[
\bar\cf_m := \text{Image of } H^*(\cf_m) \to H^*(C\u. ).
  \]
The filtration $\cf$ of $C\u.$ gives rise to the cohomology
spectral sequence $\{ E_r\}_{r\geq 1}$ converging to $H^*(C\u.
)=H^*(G/P,\bc )$.  By [Ku2, Proof of Proposition 3.2.11], for any $m\geq
0$,
  \[
E^m_1 = \bigoplus_{s+t=m} [H\u._{(s)}(\fu_P)\otimes
H\u._{(t)}(\fu^-_P)]^{\fl},
  \]
where $H\u._{(s)}(\fu_P)$ denotes the cohomology of the subcomplex
$(\wed\u._{(s)}(\fu_P))^*$ of the standard cochain complex
$\wed\u. (\fu_P)^*$ associated to the Lie algebra $\fu_P$ and
similarly for $H\u._{(t)}(\fu^-_P)$.  Moreover, by loc. cit., the
spectral sequence degenerates at the $E_1$ term, i.e.,
  \begin{equation} \label{eqn10.1}
E_1^m = E^m_{\infty}.
  \end{equation}
Further, by the definition of $P$-relative height,
  \[
[H\u._{(s)}(\fu_P) \otimes H\u._{(t)}(\fu^-_P)]^{\fl} \neq 0 \Rightarrow
s=t.
  \]
Thus,
  \begin{align*}
E^m_1 &= 0, \qquad\quad\text{unless $m$ is even and}\\
E_1^{2m} &= [H\u._{(m)}(\fu_P) \otimes H\u._{(m)}(\fu^-_P)]^{\fl} .
  \end{align*}
In particular, from (\ref{eqn10.1}) and the general properties of spectral sequences
(cf. [Ku2, Theorem E.9]), we have a canonical algebra isomorphism:
  \begin{equation} \label{eqn10.2}
\gr (\bar\cf ) \simeq \bigoplus_{m\geq 0} \bigl[ H\u._{(m)}(\fu_P) \otimes
H\u._{(m)}(\fu^-_P)\bigr]^{\fl} ,
  \end{equation}
where $\bigl[ H\u._{(m)}(\fu_P) \otimes H\u._{(m)}(\fu^-_P)\bigr]^{\fl}$
sits inside $\gr(\bar\cf )$ precisely as the homogeneous part of degree
$2m$; homogeneous parts of $\gr(\bar\cf )$ of odd degree being zero.

Finally, we claim that, for any $m\geq 0$,
  \begin{equation}  \label{eqn10.3}
\ca_m = \bar\cf_{2m} :
  \end{equation}

Following Kostant [Ko2], take the d-$\partial$ harmonic representative
 $\hat{s}^w$ in $C\u.$ for the cohomology class $[\bar{\Lambda}^P_w]$.  An
explicit expression is given by [Ko2; Theorem 4.6] together with [KK,
Theorem 3.1].  From this explicit expression, we easily see that
  \begin{equation}  \label{eqn10.4}
\ca_m \subset \bar\cf_{2m} .
  \end{equation}
Moreover, from the definition of $\ca$, for any $m\geq 0$,
  \[
\dim \frac{\ca_m}{\ca_{m+1}} = \# \bigl\{ w\in W^P : \chi_w(x_P) = \rho
+w^{-1}\rho (x_P)=m\bigr\} .
  \]
Also, by the isomorphism (\ref{eqn10.2})  and Theorem \ref{Th9.A},
  \begin{align*}
\dim \frac{\bar\cf_{2m}}{\bar\cf_{2m+1}} &= \# \bigl\{ w\in W^P: \rho-w^{-1}\rho
(x_P) = m\bigr\} \\
&= \# \bigl\{ w\in W^P: (\rho -(w_oww_o^P)^{-1}\rho )(x_P) =m\bigr\}, \\
&\qquad\qquad\text{using the involution $w \mapsto w_oww_o^P$ of $W^P$} \\
&= \bigl\{ w\in W^P: \rho +w^{-1}\rho (x_P) =m\bigr\}, \text{ since
$w_o^P$ keeps $x_P$ fixed}.
  \end{align*}
Thus,
  \[
\dim \frac{\ca_m}{\ca_{m+1}} = \dim \frac{\bar\cf_{2m}}{\bar\cf_{2m+1}} .
  \]
Further, for large enough $m_o$, $\ca_{m_o} = \bar\cf_{2m_o} =0$.  Thus, by
decreasing induction on $m$, we get (\ref{eqn10.3}) from (\ref{eqn10.4}).

Thus, combining the isomorphisms $\varphi$ and (\ref{eqn10.2}) and using
(\ref{eqn10.3}), we get the isomorphism $\phi$ as in the Theorem.  The
assertion (\ref{eqn10.0}) follows from the description of the map $\varphi$
and the identification (\ref{eqn10.2}) together with the explicit
expression of the d-$\partial$ harmonic representative $\hat{s}^w$ in
$C\u.$.  This proves the theorem.
  \end{proof}

  \begin{corollary}  \label{leviprod} The product in $(H^*(G/B), \odot_0)$
is given by
\begin{align*}
\epsilon^B_u \odot_0 \epsilon_v^B &= \epsilon_w^B,  &&\text{if $\Phi_u\cap\Phi_v =
\emptyset$ and $\Phi_w = \Phi_u\sqcup\Phi_v$}\\
&= 0, &&\text{otherwise.}
  \end{align*}
  \end{corollary}

  \begin{proof}  By the above theorem,
  \begin{align}
\phi\bigl( \epsilon_u^B\odot_0 \epsilon_v^B\bigr) &= \phi \bigl(
[\bar{\Lambda}^B_{w_ou}]\bigr) \cdot   \phi \bigl(
[\bar{\Lambda}^B_{w_ov}]\bigr)\\
\label{eqn10.5}
&= \Bigl( \frac{i}{2\pi}\Bigr)^{\ell (u)+\ell (v)} \biggl(
\prod_{\al\in\Phi_u} \ip<\rho ,\al> \cdot \prod_{\beta\in\Phi_v} \ip<\rho
,\beta>\biggr) \xi^u\, \xi^v .
  \end{align}
The right side is clearly 0 if $\Phi_u\cap\Phi_v \neq \emptyset$.  So, let us
consider the case when $\Phi_u\cap\Phi_v =\emptyset$.  In this case, two
subcases occur:
\vskip1ex

  \begin{enumerate}
\item There exists $w\in W$ such that $\Phi_w =\Phi_u\sqcup\Phi_v$.  In
particular, $\ell (w) = \ell (u)+\ell (v)$.  (Such a $w$  is necessarily
unique.)
\item There does not exist any such $w\in W$.
  \end{enumerate}

In the first case, $\xi^u\xi^v =\xi^w$ and thus the right side of equation
(\ref{eqn10.5}) is equal to $\bigl(\frac{i}{2\pi}\bigr)^{\ell (w)}\xi^w
\prod_{\al\in\Phi_w}\ip<\rho ,\al> = \phi (\epsilon_w^B)$.  Hence, in this
case, $\epsilon^B_u\odot_0\epsilon_v^B =\epsilon^B_w$.

In the second subcase, by Theorem \ref{Th9.A},
  \[
\xi^u\xi^v =0, \quad \text{ as an element of } [H^*(\fu_{B})\otimes
H^*(\fu^-_{B})]^H.
  \]
Thus, $\epsilon_u\odot_0\epsilon_v =0$.  This proves the corollary.
  \end{proof}

  \begin{remark}
a) For any $u,v\in W$, $\Phi_u\cap\Phi_v =\emptyset$ iff $\ell (uv^{-1}) =\ell
(u)+\ell (v)$.

\vskip1ex

b) A subset $S\subset R^+$ is called {\em closed under addition} if for
$\al ,\beta\in S$ such that $\al +\beta\in R^+$, we have $\al +\beta\in
S$.  A subset $S\subset R^+$ is called {\em coclosed under addition} if
$R^+\backslash S$ is closed under addition.

Now, by a result of Kostant [Ko1, Proposition 5.10], a subset $S$ of
$R^+$ is closed and coclosed under addition iff there exists $w\in W$ such
that $S=\Phi_w$.

\vskip1ex

c) For any $u,v,w\in W$, the following two conditions
are equivalent:

c$_1$) $\Phi_u\cap \Phi_v =\emptyset$ and $\Phi_w =\Phi_u\sqcup\Phi_v.$

c$_2$) $\ell (w)=\ell (u)+\ell (v)$, $\ell (uv^{-1})=\ell (u)+\ell
(v)$, $\ell (wu^{-1}) =\ell (w) -\ell (u)$ and $\ell (wv^{-1}) = \ell (w) -\ell
(v)$.

\vskip1ex

d) As a consequence of Corollary \ref{leviprod}, it is easy to get the following analogue of Chevalley formula:

For a simple reflection $s_i$,
  \begin{align*}
\epsilon^B_{s_i}\odot_0 \epsilon_v^B &= \epsilon^B_{vs_i}, &&\text{if $v\al_i$ is a
simple root}\\
 &= 0, &&\text{otherwise}.
  \end{align*}
Thus, the subalgebra of $\bigl(H^*(G/B, \Bbb Z), \odot_0\bigr)$ generated by
$H^2(G/B, \Bbb Z)$ is precisely equal to $\bigoplus\,\Bbb Z\epsilon^v$, where
the summation runs over those $v\in W$ which can be written as a product of commuting simple reflections with no simple reflection appearing more than once.
  \end{remark}
\section{Tables of the Deformed Product $\odot$ for Rank 3 Groups}\label{exemples}
We give below  the multiplication tables under the deformed product
$\odot$ (cf. equation (\ref{EE})) for $G/P$ for all the rank $3$ complex simple groups $G$ and maximal parabolic subgroups $P$. We will freely follow the convention as in [KLM] without explanation. Since we are only considering maximal parabolics, we have only one indeterminate, which we denote by ${\defpar}$. In the case of $G=$ SL$_n$, all the maximal parabolic subgroups are minuscule. Similarly, for $G=B_3$ the maximal parabolic $P_1$ is minuscule and for $G=C_3$, $P_3$ is minuscule. So, by Lemma \ref{minuscule}, the deformed product in the cohomology of the corresponding flag varieties coincides with the usual cup product, so we do not write them here. (The interested reader can find it, e.g., in [KLM].)

{\bf Example 1. $G=B_3, P=P_2:$}
\begin{center}
\begin{tabular} {|c|c|c|c|c|c|}
\hline $H^*(G/P_2)$ & $b_1$ & $b_2^{\prime}$ & $b_2^{\prime
\prime}$ & $b_3^{\prime}$ & $b_3^{\prime \prime}$ \\ \hline

$b_1$ & $b_2^{\prime} + 2 b_2^{\prime \prime}$ & $2 b_3^{\prime}$
& $b_3^{\prime} + b_3^{\prime \prime}$ & $ 2{\defpar} b_4^{\prime} +{\defpar}
b_4^{\prime \prime}$ & $ {\defpar}b_4^{\prime} + 2{\defpar} b_4^{\prime \prime}$  \\
\hline

$b_2^{\prime}$ & & $2 {\defpar}b_4^{\prime}$ & ${\defpar}b_4^{\prime} + {\defpar}b_4^{\prime
\prime}$ & $2{\defpar} b_5^{\prime} +{\defpar} b_5^{\prime \prime}$
& $ {\defpar}b_5^{\prime \prime}$ \\
\hline

$b_2^{\prime \prime}$ & & & $ {\defpar}b_4^{\prime} + {\defpar}b_4^{\prime \prime}$
& ${\defpar}b_5^{\prime} + {\defpar}b_5^{\prime \prime}$ & ${\defpar} b_5^{\prime} + {\defpar}
b_5^{\prime \prime} $  \\ \hline

$b_3^{\prime }$ & & & & $ 2 {\defpar}b_6$ & ${\defpar}b_6$ \\ \hline

$b_3^{\prime \prime}$ & & & & & $ 2 {\defpar}b_6$ \\ \hline

\end{tabular}
\end{center}

\begin{center}
\begin{tabular} {|c|c|c|c|c|c|c|}
\hline $H^*(G/P_2)$ & $b_4^{\prime}$ & $b_4^{\prime \prime}$ &
$b_5^{\prime}$ & $b_5^{\prime \prime}$ & $b_6$ & $b_7$ \\ \hline

$b_1$ & $2 b_5^{\prime}+b_5^{\prime \prime}$ & $b_5^{\prime
\prime}$ & $b_6$ & $2 b_6$ & $b_7$ & $0$ \\ \hline

$b_2^{\prime}$ & $2 b_6$ & $0$ & $b_7$ & $0$ & $0$ & $0$ \\
\hline

$b_2^{\prime \prime}$ & $b_6$ & $b_6$ & $0$ & $b_7$ & $0$ & $0$ \\
\hline

$b_3^{\prime }$ & $b_7$ & $0$ & $0$ & $0$ & $0$ & $0$\\ \hline

$b_3^{\prime \prime}$ &$0$ & $b_7$ & $0$ & $0$ & $0$ & $0$
\\ \hline
\end{tabular}
\end{center}

{\bf Example 2. $G=B_3, P=P_3:$}

\begin{center}
\begin{tabular} {|c|c|c|c|c|c|c|c|}
\hline

$H^*(G/P_3)$ & $b_1$ &  $b_2$ &  $b_3^{\prime}$ & $b_3^{\prime
\prime}$ & $b_4$ & $b_5$ &  $b_6$   \\ \hline $b_1$ & $ {\defpar}b_2$ & $
{\defpar}b_3^{\prime} + b_3^{\prime \prime}$ & $b_4$ & ${\defpar}b_4$ & ${\defpar}b_5$ &
$b_6$ & $0$ \\ \hline $b_2$ & & $2b_4$ & $b_5$ & ${\defpar}b_5$ & $b_6$ &
$0$ & $0$ \\ \hline $b_3^{\prime}$ & & & $0$ & $b_6$ & $0$ & $0$ &
$0$ \\ \hline $b_3^{\prime \prime}$ & & & & $0$ & $0$ & $0$ & $0$
\\ \hline

\end{tabular}
\end{center}

{\bf Example 3. $G=C_3, P=P_1:$}

\begin{center}
\begin{tabular} {|c|c|c|c|c|c|}
\hline

$H^*(G/P_1)$ & $a_1$ &  $a_2$ &  $a_3$ &  $a_4$ & $a_5$   \\
\hline $a_1$ & $ a_2$ & ${\defpar} a_3$ & $a_4$ &  $a_5$ &  $0$ \\ \hline
$a_2$ & & ${\defpar}a_4$ & $a_5$ &  $0$ & $0$ \\ \hline

\end{tabular}
\end{center}
{\bf Example 4. $G=C_3, P=P_2:$}
\begin{center}
\begin{tabular} {|c|c|c|c|c|c|}
\hline $H^*(G/P_2)$ & $a_1$ & $a_2^{\prime}$ & $a_2^{\prime
\prime}$ & $a_3^{\prime}$ & $a_3^{\prime \prime}$  \\ \hline

$a_1$ & $a_2^{\prime} + {\defpar}a_2^{\prime \prime}$ & ${\defpar}a_3^{\prime}$ &
$a_3^{\prime} + a_3^{\prime \prime}$ & $ 2{\defpar} a_4^{\prime} +{\defpar}
a_4^{\prime \prime}$ & $ {\defpar}a_4^{\prime} + 2{\defpar} a_4^{\prime \prime}$ \\
\hline

$a_2^{\prime}$ & & ${\defpar}^2a_4^{\prime}$ & ${\defpar}a_4^{\prime} + {\defpar}a_4^{\prime
\prime}$ & $ {\defpar}^2a_5^{\prime} + {\defpar}a_5^{\prime \prime}$ & $ {\defpar}a_5^{\prime
\prime}$  \\ \hline

$a_2^{\prime \prime}$ & & & $ 2 a_4^{\prime} + 2 a_4^{\prime
\prime}$ & ${\defpar}a_5^{\prime} + 2 a_5^{\prime \prime}$ & $ {\defpar}a_5^{\prime}
+ 2 a_5^{\prime \prime} $  \\ \hline

$a_3^{\prime }$ & & & & $ 2 {\defpar}a_6$ & ${\defpar}a_6$ \\ \hline

$a_3^{\prime \prime}$ & & & & & $ 2{\defpar} a_6$ \\ \hline

\end{tabular}
\end{center}

\begin{center}
\begin{tabular} {|c|c|c|c|c|c|c|}
\hline $H^*(G/P_2)$ & $a_4^{\prime}$ & $a_4^{\prime \prime}$ &
$a_5^{\prime}$ & $a_5^{\prime \prime}$ & $a_6$ & $a_7$ \\ \hline

$a_1$ & ${\defpar}a_5^{\prime}+a_5^{\prime \prime}$ & $a_5^{\prime \prime}$
& $a_6$ & ${\defpar}a_6$ & $a_7$ & $0$ \\ \hline

$a_2^{\prime}$ & ${\defpar}a_6$ & $0$ & $a_7$ & $0$ & $0$ & $0$ \\ \hline

$a_2^{\prime \prime}$ & $a_6$ & $a_6$ & $0$ & $a_7$ & $0$ & $0$ \\
\hline

$a_3^{\prime }$ & $a_7$ & $0$ & $0$ & $0$ & $0$ & $0$\\ \hline

$a_3^{\prime \prime}$ &  $0$ & $a_7$ & $0$ & $0$ & $0$ & $0$
\\ \hline
\end{tabular}
\end{center}

\begin{remark} We believe that our results can be generalized to the small quantum cohomology of homogenous spaces and the multiplicative eigenvalue problem. In particular, there should be an analogous deformation $\odot^q_0$ of the quantum cohomology  of homogenous spaces $G/P$ with an analogous relation to the
multiplicative eigenvalue problem ([AW],~\cite{Belkale},~\cite{Belkale4},~\cite{TW}).
\end{remark}
\bibliographystyle{plain}
\def\noopsort#1{}

  \end{document}